\UseRawInputEncoding
\documentclass[12pt]{amsart}
\usepackage[english]{babel}
\usepackage{amsfonts,amsmath,amssymb,amsthm,graphicx,afterpage,url}
\input{xy}
\xyoption{all}
\usepackage{epsfig,overpic,pgf,tikz}
\usetikzlibrary{arrows}

\usepackage{hyperref}
\hypersetup{
   colorlinks   = true, 
   urlcolor     = blue, 
   linkcolor    = blue, 
   citecolor   = red 
}

\graphicspath{{morefig/}{figures/}{figures00/}{figures02/}{figures05/}{figlms/}{figumisha/}}

\def\conv{\mathop{\fam0 conv}}

\def\t{\widetilde}
\def\R{{\mathbb R}} \def\Z{{\mathbb Z}}

\long\def\comment#1\endcomment{}

\newcommand{\jonly}[1]{}
\newcommand{\aronly}[1]{#1}

\theoremstyle{theorem}
\newtheorem{theorem}{Theorem}[section]
    \newtheorem{lemma}[theorem]{Lemma}
    
    \newtheorem{example}[theorem]{Example}
    \newtheorem{proposition}[theorem]{Proposition}
    \newtheorem{conjecture}[theorem]{Conjecture}

\newtheoremstyle{mydefinition}
  {3pt}
  {3pt}
  {\normalfont}
  {\parindent}
  {\bfseries}
  {.}
  { }
  {}

\theoremstyle{mydefinition}

\newtheorem{remark}[theorem]{Remark}

\begin{document}
 
 
\title{Realizability of hypergraphs\\ and intrinsic linking theory}

\author{A. Skopenkov}

\thanks{Homepage: \url{https://users.mccme.ru/skopenko/}.
\newline
This paper is based on the author's lectures at Moscow Institute of Physics and Technology, Independent University of Moscow, Institute of Science and Technology (Austria), and at various summer schools and math circles.
I am grateful to M. Skopenkov and A. Zimin for allowing me to use their materials, and to G. Chelnokov, I. Izmestiev, R. Karasev, 
A. Matushkin, A. Rukhovich, A. Shapovalov, M. Skopenkov, A. Sossinsky, S. Tabachnikov, O. Viro, A. Zimin, J. Zung for useful discussions.}

\date{}

\maketitle

\begin{abstract}
In this expository paper we present short simple proofs of Conway-Gordon-Sachs' theorem on intrinsic linking in three-dimensional space, as well as van Kampen-Flores' and Ummel's theorems on intrinsic intersections.
The latter are related to nonrealizability of certain hypergraphs
in four-dimensional space.
The proofs use a reduction to lower dimensions which allows to exhibit relation between these results.
We use elementary language
which allows to present the main ideas without technicalities.
Thus our exposition
is accessible to non-specialists in the area, including students who know basic three-dimensional geometry, and who are ready to learn straightforward four-dimensional generalizations.
\end{abstract}

\tableofcontents

\newpage
\hfill{\it `It's too difficult.'}

\hfill{\it `Write simply.'}

\hfill{\it `That's hardest of all.'}

\hfill{\it I. Murdoch, The Message to the Planet.}

\section{Introduction}\label{0intr}

\subsection{Impossible constructions, intrinsic intersection and intrinsic linking}

\begin{figure}[h]
\includegraphics[scale=0.7]{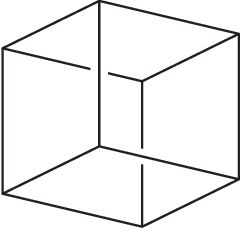} \qquad  \includegraphics[scale=0.6]{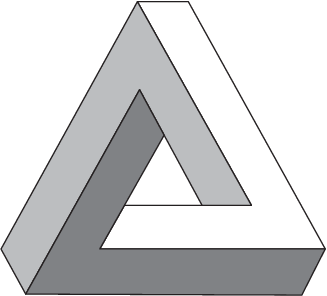} \qquad  \includegraphics[scale=0.7]{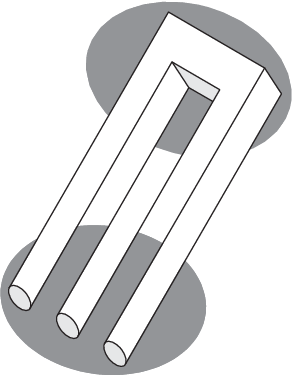} \qquad  \includegraphics[scale=0.9]{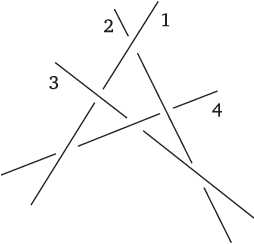}
\caption{The impossible cube, the Penrose triangle, the blivet, an impossible projection}
\label{f-imp}
\end{figure}

`Impossible constructions' like the impossible cube, the Penrose triangle, the blivet etc. are well-known, mainly due to pictures by Maurits Cornelis Escher, see Figure~\ref{f-imp}, \cite{Io}, and also \cite{Br68, CKS+, GSS+}.
The pictures do not allow the global spatial interpretation because of collision between local spatial interpretations to each other.
In geometry, topology and graph theory there are also famous basic examples of `impossible constructions'
(of which local parts are `possible').

\begin{figure}[h]\centering
\includegraphics{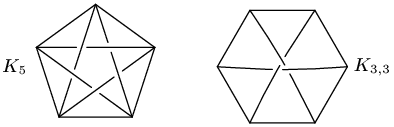}
\caption{Nonplanar graphs $K_5$ and $K_{3,3}$}
\label{k5}
\end{figure}

The following example of an `impossible construction' or an `intrinsic intersection' is already relevant to this paper. 
{\it For 5 points in the plane one cannot join each point to each other by a path so that the paths intersect only at their starting points or endpoints.}\footnote{Also, one cannot take 3 houses and 3 wells in the plane
and join each house to each well by a path so that the paths intersect only at their starting points or endpoints.
In graph-theoretic terms these assertions mean that the complete graph $K_5$ on $5$ vertices, and the complete bipartite graph $K_{3,3}$ are not planar, see  Figure~\ref{k5}.
Proposition \ref{0-ra2} below is a `linear' version of the non-planarity of $K_5$.}


\begin{proposition}[a theorem on intrinsic intersection in dimension 2]\label{0-ra2} 
For any $5$ points in the plane there are two intersecting segments joining these points, and having no common vertices.
\end{proposition}

For next results we need some notation.
We abbreviate `three-dimensional Euclidean space $\R^3$' to `3-space'.
Analogous meaning have `4-space' ($\R^4$) and $d$-space ($\R^d$).
By a \emph{triangle} we mean the part of the plane bounded by a closed polygonal line of three segments.

\begin{figure}[h]\centering
\includegraphics{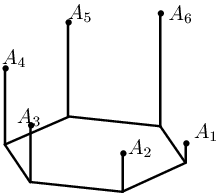} \qquad  \includegraphics{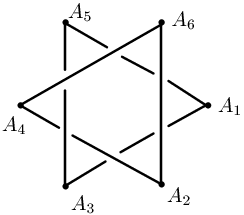} \qquad \includegraphics{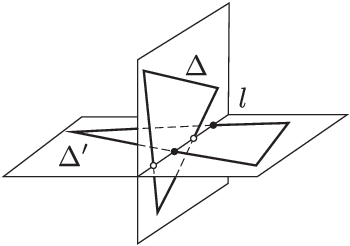}
\caption{Linked triangles}
\label{genpos}
\end{figure}

Take two triangles in 3-space no 4 of whose 6 vertices lie in one plane.
The triangles are called {\bf linked}, if the outline of the first triangle intersects
the second triangle exactly at one point.
E.g. the triangles $A_1A_3A_5$ and $A_2A_4A_6$, $\Delta$ and $\Delta'$ in Figure \ref{genpos} are linked.\footnote{The distance from the point $A_j$ to the projection plane equals $j$.
So the projection in Figure \ref{genpos}, middle, is realizable, as opposed to Figure~\ref{f-imp}, right.
\newline
The property of being linked is symmetric (this is not obvious from the definition but does have a simple proof \cite[Problem 4.1.1]{Sk}).
Other introductory results on linked triangles are presented in \cite[\S4.1 `Linked triangles']{Sk}.}


\begin{theorem}[intrinsic linking in dimension 3; Conway--Gordon--Sachs \cite{Sa81, CG83}]\label{rlt-cgslin}
If no 4 of 6 points in 3-space lie in one plane, then there are two linked triangles with vertices at these 6 points.
\end{theorem}

\begin{theorem}[intrinsic intersection in dimension 4; van Kampen-Flores \cite{vK32, Fl34}]\label{0-ne4}
From any 7 points in 4-space one can choose two disjoint triples such that the two triangles with vertices at the triples intersect.
\end{theorem}

\begin{figure}[h]
\includegraphics[width=6.5cm]{real-4.mps}
\caption{Five points in $\R^3$ (realization of the complete 3-homogeneous hypergraph on $5$ vertices)}
\label{f-full5}
\end{figure}

Analogue of theorems in dimension 2 and 4 is true for 6 points in 3-space (Proposition \ref{0-ne3}.b).
Analogues of theorems in dimension 2 and 4 are false for 4 points in the plane and for 6 points in 4-space, respectively: 
for $k=1,2$ in $\R^{2k}$ take the $2k+1$ vertices and an interior point of a $2k$-simplex, cf. Figure \ref{f-full5}.

\begin{remark}[lowering of dimension]\label{s:lowering}
Often it is convenient to reduce a planar result to a one-dimensional result (i.e., to a result in a line), and a spatial result to a planar result.
Similarly, a tempting approach to a 4-dimensional result is an analogy to, or a reduction to, a spatial result.
Some examples are given at the beginning of in \S\ref{s:interest}.
Because of such `lowering of dimension' the reader not familiar with 4-space need not be scared.
No intuition on 4-dimensional space beyond given in Remark \ref{s:intu4} is required for the proof of  Theorem \ref{0-ne4} in dimension 4. 
\end{remark}

\subsection{Why this paper might be interesting?}\label{s:interest}

We exhibit 
relation between intrinsic linking and intrinsic intersection in consecutive dimensions. 
The relation generalizes to higher dimensions (\aronly{Theorem \ref{stat-il}}\jonly{\cite[Theorem 1.6]{Sk14}};
for simplicity we mention dimensions higher than 4 only in that theorem).


Specifically, the theorem in dimension 2 (Proposition~\ref{0-ra2}) is reduced (in \S\ref{0pla}) to the following theorem in dimension 1.

\begin{proposition}[easy; a theorem on intrinsic linking in dimension 1]\label{p:liline} 
Every 4 points in a line can be colored in two red and two blue so that they alternate: 
red-blue-red-blue or blue-red-blue-red (one says `the red pair is linked with the blue pair').
\end{proposition}


Analogously, Theorem \ref{rlt-cgslin} in dimension 3 is reduced (in \S\ref{0ramcgs}) to the theorem in dimension 2 
(more precisely, to its quantitative version, Proposition \ref{0-ra2m}).
Analogously, Theorem \ref{0-ne4} in dimension 4 is reduced (in \S\ref{0ram4}) to Theorem \ref{rlt-cgslin} in dimension 3.

The results on intrinsic intersections give a natural generalization of non-planarity of graphs:
examples of {\it two-dimensional} analogues of graphs non-realizable in 3- and 4-space.
This is explained in Remark \ref{r:hyper}.

We give a simplified exposition accessible to non-specialists in the area.
We state results in terms of systems of points.
So we do not use the notion of realizability of a hypergraph
(we do mention this notion because it is an important {\it motivation}).
For understanding most of the paper it suffices to know basic geometry of 3-space, and to know or learn straightforward 4-dimensional generalizations.
We believe that the elementary description of simple applications of topological methods makes
these methods more accessible.
Comparison with other proofs is discussed in \aronly{Remark \ref{s:compar}}\jonly{\cite[Remark 3.2]{Sk14}}.

The described relation between `intrinsic intersection' and `intrinsic linking' not only gives simple proofs of classical results.
It also brings a reader to the frontline of research, notably to the solution of the generalized Menger conjecture (explained in
\S\ref{s:menger}).
The exposed intrinsic linking results are the departure point of {\it intrinsic (Ramsey) linking theory}.
\aronly{See surveys \cite{RA05, PS05, FMM+, Na20} and references therein; for higher-dimensional analogues see \cite{SS92, Sk03, KS20}.}
The exposed intrinsic intersection results are generalized to {\it non-realizability of hypergraphs}.
\jonly{See references in \cite[Remark 1.2]{Sk14}.}
\aronly{See surveys  \cite[\S4, \S5]{Sk06}, \cite[\S1]{MTW}, \cite[\S3.2]{Sk18} and references therein;
for recent results see \cite{Pa15, Sk18o, AKM, Pa21, Me22}.
For analogous problem on embedding dynamical systems see \cite{LT14} and references therein.}

The history is exposed in \aronly{Remark \ref{r:hist}}\jonly{\cite[Remark 3.6]{Sk14}}.

{\bf On the independence.}
The remarks are not formally used later and so could be omitted.
The same is true for \S\ref{s:stat4} and \S\ref{s:imrem}.
\aronly{Sections \ref{0pro}, \ref{s:imrem} and \ref{s:prod} are independent of each other,
so they could be read in any order.}
Forward references and references to other papers can be ignored for the first reading.

\subsection{Generalizations to higher dimensions and higher multiplicity}\label{s:stat4}


A subset of $\R^d$ is called {\it convex}, if for any two points from this subset the segment joining these two points is in this subset.
The {\it convex hull} of $X\subset\R^d$ is the minimal convex set that contains $X$.

\begin{theorem}\label{stat-il}
Take any $d+3$ points in $\R^d$ of which no $d+1$ points lie in one $(d-1)$-dimensional hyperplane.

For $d$ even there are two disjoint $(d+2)/2$-element subsets whose convex hulls intersect.

If $d$ is odd, then
there is an unordered pair of $(d+1)/2$-simplices with vertices at these points which is linked (i.e.,
the boundary of the first simplex intersects the convex hull of the second simplex exactly at one point).
\end{theorem}

This is Proposition \ref{0-ra2} for $d=2$, Theorem \ref{rlt-cgslin} for $d=3$, and Theorem~\ref{0-ne4} for $d=4$.
Theorem \ref{stat-il} is due to van Kampen-Flores for even $d>2$ \cite{vK32, Fl34}, and to
Segal-Spie\. z-Lov\'asz-Schrijver-Taniyama for odd $d>3$ \cite[Corollary 1.1]{LS98}, \cite{Ta00} (the index argument of \cite[\S1]{SS92} has a simple generalization to Theorem \ref{stat-il};
thus Theorem \ref{stat-il} for odd $d>3$ is implicit in \cite{SS92}).

Theorem \ref{stat-il} is proved by induction on $d$.
The base is $d=1$ and is trivial.
The inductive step is proved in \S\ref{0pro} for $d=2,3,4$; the proof for the general case is analogous.

The analogue of Theorem \ref{stat-il} for $d+2$ points

$\bullet$ does not make sense for $d$ odd because $(d+1)/2$-simplex has $(d+3)/2$ vertices;

$\bullet$ is false for $d$ even, analogously to the corresponding counterexample to Theorem~\ref{0-ne4}.

For $d$ odd there is Proposition \ref{0-ne3}.b on intrinsic intersection and its higher-dimensional analogue.
They are weaker than the corresponding Theorems \ref{rlt-cgslin}, \ref{stat-il} on intrinsic linking.
More results are presented in \cite[\S3]{RRSl}, \cite[\S4]{Sk16}.



Let us formulate a low-dimensional analogue of the above results for 3-fold intrinsic intersections.
 
\begin{theorem}[\cite{Sa91g}]\label{t:vkfr}
From any 11 points in 3-space one can choose 3 triangles having pairwise disjoint vertices but having a common point.
\end{theorem}

It is surprising that known proof of such an elementary result involves algebraic topology.
It would be interesting to obtain an elementary proof.

The analogue of Theorem \ref{t:vkfr} for 10 points is false. 
Indeed, in 3-space take the vertices of a 3-dimensional simplex and its center, see Figure \ref{f-full5}.
For every of these 5 points either take it with multiplicity two or take a close point.

For a higher-dimensional $r$-fold analogue of Theorems \ref{stat-il} on intersection and~\ref{t:vkfr} see
\cite[Theorem 1.6]{Sk16}, \cite[Conjecture 3.1.4 and the text below]{Sk18}.
For an analogue on $r$-fold linking see \cite[\S5]{AKS}.

\subsection{Some important remarks}\label{s:imrem}


\begin{remark}[Relation to hypergraphs]\label{r:hyper}
(a) Two-dimensional analogues of graphs are \emph{3-homogeneous}, or \emph{2-dimensional} \emph{hypergraphs} defined as collections of 3-element subsets of a finite set.\footnote{In topology such objects are called \emph{pure, or dimensionally homogeneous}, 2-dimensional \emph{simplicial complexes}. The term `hypergraph' is more convenient to generic mathematician or computer scientist.}
For brevity, we omit `3-homogeneous' and `2-dimensional'.
For instance, \emph{a complete hypergraph} on $k$ vertices is the collection of all 3-element subsets of a $k$-element set.
\emph{Realizability} (also called embeddability) of a hypergraph in $\R^d$ is defined similarly to the realizability of a graph in the plane: one `draws' a triangle for every three-element subset.
See Figures \ref{f-full5}, \ref{f-mnre} and~\ref{f:k5i}; on the last figures a subdivision of quadrilaterals analogous to Figure~\ref{f-mnre}, left, is not shown.
See rigorous definitions e.g. in \cite[\S3.2]{Sk18}.

Hypergraphs (and simplicial complexes) play an important role in mathematics.
One cannot imagine topology and combinatorics without them.
They are also used in computer science and in bioinformatics, see e.g. \cite{PS11}.

A `small shift' (or `general position') argument shows that every graph is realizable in 3-space.
A straightforward generalization shows that every hypergraph is realizable in 5-space.

The complete hypergraph on 6 vertices contains `the cone over $K_5$' and hence is not realizable in 3-space (Proposition \ref{0-ne3}.a).
Already in the early history of topology (1920s) mathematicians tried to construct hypergraphs non-realizable in 4-space.
Egbert van Kampen and Alejandro Flores in 1932-34 proved that the complete hypergraph on 7 vertices is not realizable in 4-space (cf. Theorem \ref{0-ne4}).
This is both an early application of {\it combinatorial topology} (nowadays called algebraic topology) and one of the first results of {\it topological combinatorics} (also an area of ongoing active research).

(b) Realizations (=embeddings) are maps without self-intersections.
For topological combinatorics and discrete geometry it is interesting to study maps whose self-intersections
are non-empty (like for embeddings), but `not too complicated'.
An important particular case is studying maps \emph{without triple intersections} and, more generally,
maps \emph{without $r$-tuple intersections}, see \S\ref{s:stat4} 
and surveys \cite{Sk16}, \cite[\S3.3]{Sk18}.

(c) We present \emph{linear} versions of the results.
\emph{PL (piecewise-linear)} and \emph{topological} realizations (=embeddings) of hypergraphs are defined and discussed e.g. in \cite[\S3.2]{Sk18}, \cite[\S6]{Sk}.
The proofs we expose are interesting because they easily generalize to the PL versions \cite{Sk03, Z13}, as opposed to the proofs of \cite{BM15, So12}.

PL versions of `quantitative' results (see \S\ref{s:quant}) imply the PL versions for \emph{almost-embeddings} (see the PL case of \cite[Theorems 1.4.1 and 3.1.6]{Sk18}).
The latter imply the \emph{topological} versions (see explanation in \cite[the paragraph after Theorem 1.4.1]{Sk18}).
\end{remark}

\begin{remark}[Comparison with other proofs]\label{s:compar}
Theorem \ref{stat-il} for $d$ even (and so its particular cases, Proposition \ref{0-ra2} and Theorem \ref{0-ne4})
has an alternative simple proof using {\it the van Kampen number}, see e.g.  \cite[\S1.4]{Sk18}, \cite[\S1.4, \S6]{Sk}.
(That proof works for quantitative, PL, and topological versions, see Remark \ref{r:hyper}.c.)
That proof and the proof sketched in this paper, are presumably the simplest known proofs (`proofs from the Book').

Usually Theorem \ref{stat-il} for $d$ even (more precisely, the topological version of Theorems \ref{stat-il} and \ref{0-ne4j}) is proved using the Borsuk-Ulam theorem \cite[\S8]{Sk20}, \cite[\S5]{Ma03}.
As opposed to this paper (and to the alternative proof using the van Kampen number), this requires some knowledge of algebraic topology.
This knowledge does not make things simpler: known proofs of the Borsuk-Ulam theorem (see \cite{Ma03} and the references therein) are not easier than the above-discussed direct proofs of Theorem \ref{stat-il} for $d$ even.
(The Borsuk-Ulam theorem is proved using {\it the degree} analogously to the direct proof of Theorem \ref{stat-il} for $d$ even using {\it the van Kampen number}.)

The paper \cite{BM15} presents short algebraic proof of Theorem \ref{stat-il} (and so of its particular cases).
That proof is in the spirit of the algebraic proof of the {\it Radon theorem}\footnote{See e.g. \cite[\S1]{Sk16} for the statement of the Radon theorem.
See \cite[\S4]{Sk16} for relations between the Radon theorem and Theorems \ref{rlt-cgslin}, \ref{0-ne4}, \ref{stat-il}.
The proof of \cite{BM15} is presumably a direct (i.e., without use of the Gale transform) version of the proof of
\cite[Theorem 5]{So12} for $k=1$, which is Theorem \ref{stat-il} for $d$ even.}.
\end{remark}

\begin{remark}[history]\label{r:hist}
General `lowering of dimension' or `the link of a vertex' ideas are simple and well-known (see Remark \ref{s:lowering}).
For proofs of the Radon theorem based on this idea see \cite{Pe72, Ko18, RRS}.
For an application in computer science see \cite[proof of 2.3.i]{DE94}.
Also well-known is relation between linking and intersection.\footnote{E.g. the linking number of two disjoint closed polygonal lines in 3-dimensional sphere $\partial D^4$ equals to the algebraic intersection number of two \emph{general position} 2-dimensional disks in 4-dimensional ball $D^4$ spanning the two polygonal lines.
For a certain inductive argument involving assertion on linking in odd dimensions and assertion on intersection in even dimensions see \cite[Whitney Lemma 5.12 and Theorem 5.16]{RS72}.}
An elaboration of this idea to a relation between intrinsic linking and non-realizability
is non-trivial (cf. the difference between Proposition \ref{0-ne3}.a and Theorem \ref{rlt-cgslin}).
Proofs that discover and use that relation seem to have not been published

$\bullet$ before \cite{RST91, RST93, Sh03, RSSZ, Z13}, for a proof of the Conway--Gordon--Sachs Theorem \ref{rlt-cgslin} by reducing intrinsic linking to intrinsic intersection in lower dimension,


$\bullet$ before \cite[Example 2, Lemmas 2 and 1']{Sk03}, \cite{RSSZ}, for proofs of Theorem \ref{0-ne4} and of the Menger conjecture (see \S\ref{s:menger}) by reducing intrinsic intersection  to intrinsic linking in lower dimension.
\end{remark}

\section{Proofs}\label{0pro}

By {\it $k$ points in $\R^d$} (in this paper mostly $d\le4$) we mean $k$ pairwise distinct points.
Denote $[n]:=\{1,2,\ldots,n\}$.

\subsection{Intersections in the plane: proof of Proposition \ref{0-ra2}}\label{0pla}

Proposition \ref{0-ra2} is easily proved by analyzing the convex hull of the 5 points.
In order to illustrate the `lowering of dimension' idea (see Remark \ref{s:lowering}) in the simplest situation,
we also deduce Proposition \ref{0-ra2} from Proposition \ref{p:liline}.  


\begin{figure}[h]\centering
\includegraphics{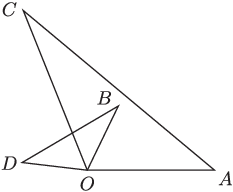} \aronly{\qquad \qquad\includegraphics[scale=0.6]{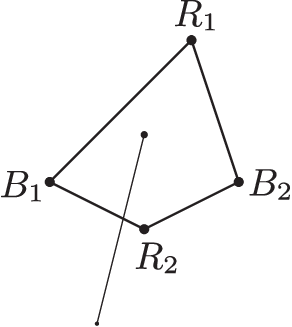}}
\caption{\aronly{Left:} To the proof of Proposition~\ref{0-ra2}\aronly{. Right: To Proposition~\ref{0-k33}.b.}}
\label{f:plane}
\end{figure}

\begin{proof}[Proof of Proposition \ref{0-ra2}]
There is a line $l$ such that one of the given points, say $O$, lies on one side of $l$, and the other points, say   $A,B,C,D$, lie on the other side of $l$.
If for some two points $X,Y\in\{A,B,C,D\}$ the point $X$ belongs to the segment $OY$, then we are done.
Otherwise we can assume that the points $A,B,C,D$ are seen from $O$ in this order, see Figure \ref{f:plane}.
Then the required assertion follows by Lemma \ref{0-eveng} below.
\end{proof}


\begin{lemma}[lowering of dimension; see Figure~\ref{f:plane}\aronly{, left}]\label{0-eveng}
Two triangles in the plane have a common vertex.
A line $l$ splits this point from the bases of the triangles.
The intersections of $l$ with the triangles alternate along $l$.
Then some two sides of the triangles intersect but do not have common vertices.
\end{lemma}

This lemma easily follows from Proposition \ref{p:liline}.
The lemma is explicitly stated in order to conveniently use it (here and in \S\ref{0ram4}), and to
illustrate its generalization to higher-dimensions (\aronly{Lemma \ref{0-even4g}}\jonly{\cite[Lemma ???]{Sk14}}).

\jonly{See more results in \cite[\S2.1]{Sk14}.}

\aronly{

The following propositions are proved analogously to Proposition \ref{0-ra2}.
They are used for some 3-dimensional results (Propositions \ref{p:k4k4}.ab and Theorem \ref{rlt-sachs}).

\begin{proposition}\label{0-k33}
(a) (Figure~\ref{k5}, right) Two triples of points are given in the plane.
Then there exist two intersecting segments without common vertices and such that each segment joins the points from distinct triples.

(b) (Figure~\ref{f:plane}, right) Four red and two blue points $B_1, B_2$ are given in the plane.
Suppose that any two segments joining points of different colors either are disjoint or intersect at their common vertex.
Then there are two red points $R_1,R_2$ such that the quadrilateral $R_1B_1R_2B_2$ does not have self-intersections, and the remaining two red points lie in different sides w.r.t. the quadrilateral.
(`In different sides' means that a \emph{general position} polygonal line joining the remaining two red points intersects the outline of the quadrilateral at an odd number of points.)
\end{proposition}

}

\subsection{Weaker versions of Theorem \ref{rlt-cgslin}}\label{s:weaker}

First we illustrate the `lowering of the dimension' idea (see Remark \ref{s:lowering}) 
by proving the following weaker versions of Theorem \ref{rlt-cgslin}.

\begin{proposition}\label{0-ne3} From any 6 points in 3-space one can choose

(a) 5 points $O,A,B,C,D$ such that the triangles $OAB$ and $OCD$ have a common point other than $O$.

(b) disjoint pair and triple such that the segment joining points of the pair intersects the triangle spanned by the triple.
\end{proposition}


\begin{figure}[h]\centering
\begin{tikzpicture}[line cap=round,line join=round,>=triangle 45,x=0.32cm,y=0.32cm]
\clip(-4.93,-7.56) rectangle (15.43,4.12);
\fill[fill=black,fill opacity=0.1] (-3,-5) -- (0,1) -- (13,1) -- (10,-5) -- cycle;
\draw (4,3)-- (1,0);
\draw (4,3)-- (-1,-4);
\draw (4,3)-- (3,-2);
\draw (4,3)-- (7,-2);
\draw (4,3)-- (11,0);
\draw (2.4,-5)-- (2,-7);
\draw [dash pattern=on 2pt off 2pt] (3,-2)-- (2.4,-5);
\draw [dash pattern=on 2pt off 2pt] (7,-2)-- (8.76,-5);
\draw (8.76,-5)-- (9.85,-6.86);
\draw [dash pattern=on 2pt off 2pt] (1,0)-- (-2,-3);
\draw (-2,-3)-- (-3,-4);
\draw (-1.71,-5)-- (-2.4,-5.96);
\draw (12.24,-0.53)-- (13.38,-1.02);
\draw [dash pattern=on 2pt off 2pt] (-1.71,-5)-- (-1,-4);
\draw [dash pattern=on 2pt off 2pt] (11,0)-- (12.24,-0.53);
\draw (-3,-5)-- (0,1);
\draw (13,1)-- (10,-5);
\draw (10,-5)-- (-3,-5);
\draw (0,1)-- (1.66,1);
\draw (3,1)-- (3.34,1);
\draw (4,1)-- (5,1);
\draw (6,1)-- (8,1);
\draw (9.32,1)-- (13,1);
\fill [color=black] (1,0) circle (1.5pt);
\fill [color=black] (-1,-4) circle (1.5pt);
\fill [color=black] (7,-2) circle (1.5pt);
\fill [color=black] (11,0) circle (1.5pt);
\fill [color=black] (3,-2) circle (1.5pt);
\fill [color=black] (4,3) circle (1.5pt);
\draw[color=black] (2.83,3.24) node {$O$};
\fill [color=black] (2,-7) circle (1.5pt);
\draw[color=black] (3.27,-6.36) node {$A_5$};
\fill [color=black] (9.85,-6.86) circle (1.5pt);
\draw[color=black] (8.69,-6.32) node {$A_4$};
\fill [color=black] (-3,-4) circle (1.5pt);
\draw[color=black] (-3.37,-2.64) node {$A_2$};
\fill [color=black] (-2.4,-5.96) circle (1.5pt);
\draw[color=black] (-1.09,-5.8) node {$A_3$};
\fill [color=black] (13.38,-1.02) circle (1.5pt);
\draw[color=black] (13.05,-2.14) node {$A_1$};
\end{tikzpicture} 
\caption{To the proofs of Proposition \ref{0-ne3}.a and Theorem~\ref{rlt-cgslinm}.
A plane in $\R^3$ intersects the segments $OA_1,\ldots,OA_5$ by points $A_1',\ldots,A_5'$.}
\label{f:rlt-cgslin}
\end{figure}

\begin{proof}[Proof of (a)]
There is a plane $\alpha$ such that one of the given points, say $O$, lies on one side of $\alpha$, and the other $5$ points lie on the other side of $\alpha$ (Figure \ref{f:rlt-cgslin}).
Consider the intersection of $\alpha$ with the union of triangles $OAB$ for all pairs $A,B$ of given points.
Now part (a) follows by Proposition~\ref{0-ra2}.
\end{proof}

Part (b) follows from (a) (and vice versa).


Figure \ref{f-full5} shows that the analogues of (a,b) for 5 points are both false.

\subsection{`Quantitative' versions}\label{s:quant}

We prove the following stronger `quantitative' (i.e., algebraic modulo 2) version of the results from \S\ref{0intr}.
(The deduction of Theorem \ref{rlt-cgslin} from Proposition~\ref{0-ra2}, not from its quantitative version, has a technical detail
which is hard to generalize to higher dimensions.)

{\renewcommand{\thetheorem}{\ref{p:liline}$'$}\addtocounter{theorem}{-1}
\begin{proposition}[obvious]\label{p:liline'}
Any 4 points in the line have a unique unordered splitting into two linked pairs.
\end{proposition}
}

{\renewcommand{\thetheorem}{\ref{0-ra2}$'$}\addtocounter{theorem}{-1}
\begin{proposition}\label{0-ra2m}
No 3 of given 5 points in the plane lie in one line.
Consider all segments joining these 5 points.
Then the number of those intersection points of such segments, that are distinct from the vertices is odd.
\end{proposition}
}

This is easily proved by analyzing the convex hull of the 5 points, or is implied by Proposition \ref{p:liline'} and the following lemma.

{\renewcommand{\thetheorem}{\ref{0-eveng}$'$}\addtocounter{theorem}{-1}
\begin{lemma}\label{0-evengq}
(a) Two triangles in the plane have a common vertex $O$.
No 3 of their 5 vertices lie in one line.
A line $l$ splits $O$ from the bases of the triangles.
Then the outlines of the triangles intersect at an even number of points if and only if the intersections of $l$ with the outlines alternate along $l$ (i.e., if the intersection of one triangle with the outline of the other contains exactly one segment with vertex $O$).

(b) No $3$ of $5$ points $O,A_1,A_2,A_3,A_4$ in the plane lie in one line.
Then the number from Proposition~\ref{0-ra2m} equals to the sum of the numbers of intersection points of the interiors of sides of triangles $OPQ$ and $ORS$, over all unordered splittings of points $A_1,A_2,A_3,A_4$ into two unordered pairs $P,Q$ and $R,S$.

(c) \cite[Proposition 7.5.a]{DGN+} Denote by $X={{[5]\choose2}\choose2}$ the set of all unordered pairs of $2$-element subsets of $[5]$.
    For any of the $3$ non-ordered partitions $\sigma\sqcup\tau=[4]$ into $2$-element sets denote
    \[
        T_{\{\sigma,\tau\}}:=\bigl\{\{\alpha,\beta\}\in X\ :\  \alpha\subset\sigma\sqcup\{5\},\ \beta\subset\tau\sqcup\{5\}\bigr\}.
    \]
Then a pair $\{\alpha,\beta\}\in X$ is contained in an odd number of sets $T_{\{\sigma,\tau\}}$ if and only if
$\alpha\cap\beta=\emptyset$.\footnote{In other words, the sum modulo $2$ (i.~e., the symmetric difference) of the sets $T_{\{\sigma,\tau\}}$ over all the $3$ such partitions $\sigma\sqcup\tau=[4]$ equals $\bigl\{\{\alpha,\beta\}\in X\  :\ \alpha\cap\beta=\emptyset\bigr\}$.}
\end{lemma}
}

Part (a) is analogous to lemma \ref{0-eveng} (and so is trivial).
Part (b) follows from (c).
A simple proof of (c) is left to a reader.

\begin{remark}\label{r:stronger}
Proposition \ref{0-ra2m} is indeed stronger than Proposition \ref{0-ra2} because it suffices to prove
Proposition \ref{0-ra2} under the assumption that no 3 of the 5 points lie in one line,

(a) first, since otherwise Proposition \ref{0-ra2} is obvious:
if points $A,B,C$ among given 5 points lie in one line, $B$ between $A$ and $C$, and $D$ is any other given point, then segments $AC$ and $BD$ intersect.

(b) second, since we can make a small shift so that no 3 of the 5 shifted points lie in one line, and no intersection points of segments with disjoint vertices are added.

Analogously to (b), Theorems \ref{rlt-cgslinm} and \ref{0-ne4m} below are stronger than Theorems \ref{rlt-cgslin} and \ref{0-ne4}.
\end{remark}

Triangles differing only by permutation of vertices are considered to be the same.

{\renewcommand{\thetheorem}{\ref{rlt-cgslin}$'$}\addtocounter{theorem}{-1}
\begin{theorem}[\cite{Sa81, CG83}]\label{rlt-cgslinm}
No 4 of given 6 points in 3-space lie in one plane.
Then the number of linked unordered pairs of triangles with vertices at these 6 points is odd.
\end{theorem}
}


\aronly{

{\renewcommand{\thetheorem}{\ref{0-ne4}$'$}\addtocounter{theorem}{-1}
\begin{theorem}[\cite{vK32, Fl34}]\label{0-ne4m}
No 5 of given 7 points in 4-space lie in one 3-dimensional hyperplane. 
Consider all triangles with vertices at these 7 points. 
Then the number of those intersection points of such triangles, that are not contained in the union of edges is odd.
\end{theorem}
}

Theorem \ref{stat-il} has an analogous quantitative version.

Another `quantitative versions' are presented in \S\ref{s:unl}.
For counterexamples to `integer versions' see \cite[Proposition 1.2 and Theorem 1.4]{KS20}.}

\jonly{Theorems \ref{0-ne4} and \ref{stat-il} have analogous quantitative version.}

\subsection{Linking in 3-space: proof of  Theorem \ref{rlt-cgslinm}}\label{0ramcgs}

\aronly{\begin{proof}[Proof of  Theorem \ref{rlt-cgslinm}]}
There is a plane $\alpha$ such that one of the given points, say $O$, lies on one side of $\alpha$, and 
the other points, say $A_1,\ldots,A_5$, lie on the other side of $\alpha$ (Figure \ref{f:rlt-cgslin}).
Take the intersection points of $\alpha$ and segments $OA_1,\ldots,OA_5$.
Since no 4 of the 6 points lie in one plane, no 3 of the taken 5 points lie in one line.
So looking from point $O$ in the plane $\alpha$ we obtain a picture analogous to Figure \ref{genpos}, middle.

In 3-space \textit{a segment $MN$ is below a segment $q$ (looking from point $O$)}, if $q$ intersects the interior of the triangle $OMN$.
The triangles $OA_1A_2$ and $A_3A_4A_5$ are linked if and only if $A_1A_2$ is below \emph{exactly one side} of the triangle $A_3A_4A_5$.
A segment cannot intersect a triangle by more than two points. 
So we can replace `exactly one side' by `an odd number of sides'.

Then the following numbers have the same parity:

$\bullet$ the number of linked unordered pairs of triangles formed by given 6 points;

$\bullet$ the number of segments $A_iA_j$ that are below an odd number of sides of their `complementary' triangles $A_kA_lA_m$, $\{i,j,k,l,m\}=\{1,2,3,4,5\}$;

$\bullet$ the number of `undercrossings', i.e., of ordered pairs $(A_iA_j,A_kA_l)$ of segments in which the first segment is below the second one;

$\bullet$ the number of intersection points of interiors of segments whose vertices are the 5 taken points in $\alpha$.

By Proposition~\ref{0-ra2m} the latter number is odd.
Hence the first number is also odd.
\aronly{\end{proof}}



\aronly{

The following version of Theorem \ref{rlt-cgslin} is analogously reduced to Proposition \ref{0-k33}.b \cite{Z13}.
This version is used for some 4-dimensional result (Theorem \ref{0-ne4pr}) in \S\ref{s:mnre4}.

\begin{figure}[h]\centering
\includegraphics[scale=0.8]{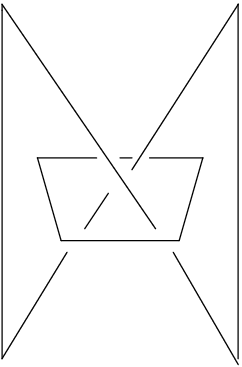}
\caption{Whitehead link formed by space quadrilaterals}
\label{f:whit}
\end{figure}

In 3-space take two quadrilaterals (i.e., closed quadrangular polygonal lines) $ABCD$ and $A'B'C'D'$,
no 4 whose 8 vertices lie in one plane.
The quadrilaterals are called {\it linked modulo 2} if the number of intersection points of the quadrilateral $ABCD$ with the union of the triangles $A'B'C'$ and $A'D'C'$ is odd.
(As opposed to triangles, there are space quadrilaterals {\it linked} but not linked modulo 2,
see Figure \ref{f:whit}.
Cf. definition at the beginning of \S\ref{s:mnre4}.)


\begin{theorem}[\cite{Sa81}]\label{rlt-sachs}
There are 4 points red and 4 blue points in 3-space.
No 4 of these 8 points lie in one plane.
Then there are two linked modulo 2 space quadrilaterals consisting of segments joining points of different colors.
\end{theorem}

}

\subsection{Intersection in 4-space: proof of Theorem \ref{0-ne4}}\label{0ram4}


\begin{remark}[some intuition on 4-space]\label{s:intu4}
(a) `Typical' intersection

$\bullet$ of two segments in the plane is either empty set or a point (here `typical' means that no 3 points among the vertices of segments lie in one line);

$\bullet$ of a segment and a triangle in 3-space is either empty set or a point;


$\bullet$ of two triangles in 4-space is either empty set or a point.

(b) For each two points

$\bullet$ of the plane distinct from a point $A$ in the plane there exists a polygonal line joining these points and not passing through $A$.

$\bullet$ of 3-space not belonging to a line $l$ in 3-space there exists a polygonal line joining these points and disjoint with $l$;


$\bullet$ in 4-space not belonging to a 2-dimensional plane $\alpha$ in 4-space
there exists a polygonal line joining these points and disjoint with $\alpha$.
\end{remark}


More intuition on 4-space is not required here, but can be developed by studying e.g. \cite[\S1]{RRSl},
\cite[\S5.1 `How to work with four-dimensional space?']{Sk}.

\begin{proof}[Proof of  Theorem \ref{0-ne4}]
There is a 3-dimensional hyperplane $\alpha$ such that one of the given points, say $O$, lies on one side of $\alpha$, and the other 6 points lie on the other side of $\alpha$ (Figure \ref{f:4space}).
Take the 6 intersection points of $\alpha$ with the segments joining $O$ to the other 6 points.
We may assume that no 5 of the given 7 points lie in one 3-dimensional hyperplane (analogously to Remark \ref{r:stronger}.b).
Hence no 4 of the 6 points in $\alpha$ lie in one plane.
Then by Theorem~\ref{rlt-cgslin} there are two linked triangles with vertices at the taken 6 points.
So we are done by Lemma \ref{0-even4}.
\end{proof}

\begin{figure}[h]\centering
\begin{tikzpicture}[line cap=round,line join=round,>=triangle 45,x=0.32cm,y=0.32cm]
\clip(-4.93,-7.56) rectangle (15.43,4.12);
\fill[fill=black,fill opacity=0.1] (-3,-5) -- (0,1) -- (13,1) -- (10,-5) -- cycle;
\draw (1,0)-- (-1,-4);
\draw (-1,-4)-- (7,-2);
\draw (9,-4)-- (11,0);
\draw (3,-2)-- (4.51,-2.5);
\draw (9,-4)-- (5.01,-2.67);
\draw (7,-2)-- (5.54,-1.51);
\draw (4,3)-- (1,0);
\draw (4,3)-- (-1,-4);
\draw (4,3)-- (3,-2);
\draw (4,3)-- (7,-2);
\draw (4,3)-- (9,-4);
\draw (4,3)-- (11,0);
\draw (2.4,-5)-- (2,-7);
\draw [dash pattern=on 2pt off 2pt] (3,-2)-- (2.4,-5);
\draw [dash pattern=on 2pt off 2pt] (7,-2)-- (8.76,-5);
\draw (8.76,-5)-- (9.85,-6.86);
\draw [dash pattern=on 2pt off 2pt] (1,0)-- (-2,-3);
\draw (-2,-3)-- (-3,-4);
\draw (-1.71,-5)-- (-2.4,-5.96);
\draw (9.71,-5)-- (10.74,-6.43);
\draw (12.24,-0.53)-- (13.38,-1.02);
\draw [dash pattern=on 2pt off 2pt] (-1.71,-5)-- (-1,-4);
\draw [dash pattern=on 2pt off 2pt] (9,-4)-- (9.71,-5);
\draw [dash pattern=on 2pt off 2pt] (11,0)-- (12.24,-0.53);
\draw (-3,-5)-- (0,1);
\draw (13,1)-- (10,-5);
\draw (10,-5)-- (-3,-5);
\draw (0,1)-- (1.66,1);
\draw (3,1)-- (3.34,1);
\draw (4,1)-- (5,1);
\draw (6,1)-- (8,1);
\draw (9.32,1)-- (13,1);
\draw (3,-2)-- (6.25,-1.19);
\draw (7.22,-0.94)-- (11,0);
\draw (1,0)-- (1.5,-0.17);
\draw (1.89,-0.3)-- (2.85,-0.62);
\draw (3.57,-0.86)-- (5.02,-1.34);
\fill [color=black] (1,0) circle (1.5pt);
\fill [color=black] (-1,-4) circle (1.5pt);
\fill [color=black] (7,-2) circle (1.5pt);
\fill [color=black] (9,-4) circle (1.5pt);
\fill [color=black] (11,0) circle (1.5pt);
\fill [color=black] (3,-2) circle (1.5pt);
\fill [color=black] (4,3) circle (1.5pt);
\draw[color=black] (2.83,3.24) node {$O$};
\fill [color=black] (2,-7) circle (1.5pt);
\draw[color=black] (3.27,-6.36) node {$A_5$};
\fill [color=black] (9.85,-6.86) circle (1.5pt);
\draw[color=black] (8.69,-6.32) node {$A_4$};
\fill [color=black] (-3,-4) circle (1.5pt);
\draw[color=black] (-3.37,-2.64) node {$A_2$};
\fill [color=black] (-2.4,-5.96) circle (1.5pt);
\draw[color=black] (-1.09,-5.8) node {$A_3$};
\fill [color=black] (10.74,-6.43) circle (1.5pt);
\draw[color=black] (11.63,-5.48) node {$A_6$};
\fill [color=black] (13.38,-1.02) circle (1.5pt);
\draw[color=black] (13.05,-2.14) node {$A_1$};
\end{tikzpicture} 
\caption{To the proof of  Theorem \ref{0-ne4};
the hyperplane $\alpha$ in 4-space is shown as a plane in 3-space}
\label{f:4space}
\end{figure}

\begin{lemma}[lowering of the dimension; Figure \ref{f:4space}]\label{0-even4}
Two tetrahedra $\tau$ and $\tau'$ in 4-space have a common vertex.
A 3-dimensional hyperplane $\alpha$ splits this vertex from the bases of the tetrahedra.
The outlines of the triangles $\alpha\cap\tau$ and $\alpha\cap\tau'$ are disjoint and linked in $\alpha$.
Then some two faces of the tetrahedra intersect but do not have common vertices.
\end{lemma}

This lemma is not as obvious as its low-dimensional analogues (Lemma \ref{0-eveng} and analogous result for a triangle and a tetrahedron in 3-space)
because {\it the surface of a tetrahedron in 4-space does not split 4-space} (cf. Remark \ref{s:intu4}.b).


\begin{proof}[Proof of Lemma \ref{0-even4}]
Denote by $\gamma$ the intersection plane of the 3-dimensional hyperplanes spanned by the tetrahedra.
Then $\alpha\cap\gamma$ is the intersection line of the planes of the linked triangles $\alpha\cap\tau$ and $\alpha\cap\tau'$.
Hence $\tau\cap\gamma$ and $\tau'\cap\gamma$ are triangles with a common vertex $O$, which is a common vertex of the tetrahedra (Figure \ref{f:plane}\aronly{, left}).
Since $\alpha$ splits $O$ from the bases of the tetrahedra, we have that $\alpha\cap\gamma$ splits $O$ from the bases of the triangles.
Since the triangles $\alpha\cap\tau$ and $\alpha\cap\tau'$ are linked, the intersection points of the line $\alpha\cap\gamma$ and the outlines of the triangles $\tau\cap\gamma$ and $\tau'\cap\gamma$ alternate along the line \cite[Problem 4.1.1]{Sk}  
(cf. Figure \ref{genpos}, right).
Hence by Lemma \ref{0-eveng} two sides of the triangles $\tau\cap\gamma$ and $\tau'\cap\gamma$ intersect but do not have common vertices.
At most one of these sides contains $O$.
Hence the two sides are contained in two faces of the tetrahedra intersect but do not have common vertices.
\end{proof}


Theorem \ref{0-ne4m} follows from Theorem~\ref{rlt-cgslinm} and the following Lemma \ref{relation}.

{\renewcommand{\thetheorem}{\ref{0-even4}$'$}\addtocounter{theorem}{-1}
\begin{lemma}\label{relation}
(a) Two tetrahedra $\tau$ and $\tau'$ in 4-space have a common vertex $O$.
No 5 of their 7 vertices lie in one 3-dimensional hyperplane.
A 3-dimensional hyperplane $\alpha$ splits $O$ from the bases of the tetrahedra.
Then the surfaces of the tetrahedra $\tau,\tau'$ intersect at an even number of points if and only if the
triangles $\alpha\cap\tau$ and $\alpha\cap\tau'$ are linked in $\alpha$ (i.e., if the intersection of one tetrahedron with the surface of the other contains exactly one segment with the end $O$).

(b) No $5$ of $7$ points $O,A_1,\ldots,A_6\in\R^4$ lie in one 3-dimensional hyperplane.
Then the number from Theorem~\ref{0-ne4m} equals to the sum of the numbers of intersection points of the interiors of faces of tetrahedra $O\Delta$ and $O\Delta'$, over all unordered splittings of points $A_1,\ldots,A_6$ into two unordered triples $\Delta$ and $\Delta'$.
(The \emph{interior} of a triangle is its complement to the outline.)

(c) Denote by $X={{[7]\choose3}\choose2}$ the set of all unordered pairs of $3$-element subsets of $[7]$.
    For any of the $10$ non-ordered partitions $\sigma\sqcup\tau=[6]$ into $3$-element sets denote
    \[
        T_{\{\sigma,\tau\}}:=\bigl\{\{\alpha,\beta\}\in X\ :\  \alpha\subset\sigma\sqcup\{7\},\ \beta\subset\tau\sqcup\{7\}\bigr\}.
    \]
Then a pair $\{\alpha,\beta\}\in X$ is contained in an odd number of sets $T_{\{\sigma,\tau\}}$ if and only if
$\alpha\cap\beta=\emptyset$.\footnote{In other words, the sum modulo $2$ of the sets $T_{\{\sigma,\tau\}}$ over all such partitions $\sigma\sqcup\tau=[6]$ equals
\linebreak
$\bigl\{\{\alpha,\beta\}\in X\  :\ \alpha\cap\beta=\emptyset\bigr\}$.}
\end{lemma}
}

Part (a) is analogous to lemma \ref{0-even4}.
Part (b) follows from (c).
A simple proof of (c) is left to a reader.

The following higher-dimensional version of Proposition \ref{0-k33}.a is related (analogously to Theorem \ref{0-ne4}) to some 3-dimensional intrinsic linking result \cite[Remark 2.5]{DS22}.

\begin{theorem}[\cite{vK32, Fl34}]\label{0-ne4j} Three triples of points in 4-space are given.
Then there exist two intersecting triangles without common vertices
such that the vertices of each triangle belong to distinct triples.
\end{theorem}

\subsection{Unlinking properties}\label{s:unl}


Here we present quantitative versions asserting that the number of intersections (or linkings) is even,
cf. \S\ref{s:quant}.

\begin{proposition}\label{r:lint}
(2) There are 5 points in the plane such that no 3 of them lie in a line, and every segment joining two of the points intersects the outline of the triangle formed by the remaining three points at an even number of points.
(I.e., every pair of points is `unlinked' with the triangle formed by the remaining three points.)

(2') No 3 of 5 given points in the plane lie in a line.
Then the number of those segments joining two of the points that intersect at exactly one point the outline of the triangle formed by the remaining three points, is even.
\end{proposition}

Proofs are easy and are left to the reader.


In 3-space instead of unlinking properties \ref{r:lint}.2,2' there is a linking property (Theorem \ref{rlt-cgslin}) and the following unlinking properties.

\begin{proposition}\label{r:lint3}
(3) There are 6 points in 3-space such that no 4 of them lie in a plane, and every segment joining two of them intersects the surface of the tetrahedron formed by the remaining four points at an even number of points.
(I.e., every pair of points is `unlinked' with the tetrahedron formed by the remaining four points.)


(3') No 4 of given 6 points in 3-space lie in one plane.
Then the number of those segments joining two of the points that intersect at exactly one point the surface of the tetrahedron formed by the remaining four points, is even.
\end{proposition}

\begin{proof}[Sketch of a proof]
(3) Take points close to the vertices of regular octahedron, points close to the vertices of a triangular prism, or points on the moment curve.


(3') Any two triangles spanned by two disjoint triples of given points either are disjoint or intersect by a segment (non-degenerate to a point).
There is an even number of ends of such segments.
The ends of such segments are exactly intersection points of segments joining pairs of points, and surfaces of `complementary' tetrahedra.

Alternatively, we have
$$\sum_{\{A,B\}} |AB\cap\partial T_{AB}|_2 = \sum_{\{A,B\}} (|A\cap T_{AB}|_2+|B\cap T_{AB}|_2) =
\sum_{A}\sum_{B\ne A} |A\cap T_{AB}|_2=0.$$
Here

$\bullet$ $|S|_2$ is the parity of the number of elements in a finite set $S$,

$\bullet$ $T_{AB}$ is the tetrahedron formed by the four given points distinct from $A,B$, and

$\bullet$ the last equality holds because for every $A$ the set  $\{T_{AB}\}_{B\ne A}$ is a \emph{3-cycle}, cf. \cite[Remark 1.3.5b]{Sk18}.
\end{proof}

Propositions \ref{0-ra2}, \ref{0-ra2m}, \ref{0-ne3}.b and \ref{r:lint3}.3'
show that under transition from dimension 2 to dimension 3 the property of the existence of intersection is preserved, while the parity of the number of intersections change.
The 3-dimensional versions of Propositions \ref{0-ra2}, \ref{0-ra2m} have a stronger form: Theorems \ref{rlt-cgslin} and \ref{rlt-cgslinm}.


\begin{proposition}\label{r:lint4}
(4-3) There are 7 points in 4-space such that no 5 of them lie in a  3-dimensional hyperplane, and every triangle formed by 3 of them intersects the surface of the tetrahedron formed by the 4 remaining points at an even number of points.
(I.e., every triangle formed by three of the points is `unlinked' with the tetrahedron formed by the remaining four points.)

(4'-3) No 5 of 7 given points in 4-space lie in a 3-dimensional hyperplane.
Then the number of those triangles spanned by three of the points that intersect at exactly one point the surface of the tetrahedron formed by the remaining four points, is even.


(4-2) (conjecture) There are 7 points in 4-space such that no 5 of them lie in a 3-dimensional hyperplane, and every segment joining two of them intersects the 3-dimensional surface of the 4-simplex formed by the remaining five points at an even number of points.
(I.e., every pair of points is `unlinked' with the 4-simplex formed by the remaining five points.)


(4'-2) No 5 of 7 given points in 4-space lie in a 3-dimensional hyperplane.
Then the number of those segments joining two of the points that intersect at exactly one point the 3-dimensional surface of the 4-dimensional simplex formed by the remaining five points, is even.
\end{proposition}

\begin{proof}[Sketch of a proof]
(4-3) Take points on the moment curve. See details in \cite{St24}.

(4'-3) Analogously to Proposition \ref{r:lint3}.3'.


(4-2) Perhaps one can take points on the moment curve.

(4'-2) Analogously to the alternative proof of Proposition \ref{r:lint3}.3'.
\end{proof}

\begin{conjecture}\label{r:lintd}
(d-k) For any $d\ne 2k-1$ there are $d+3$ points in $\R^d$, of which no $d+1$ lie in one $(d-1)$-hyperplane, and such that any $k$-simplex spanned by $k+1$ of them, intersects the surface of $(d+1-k)$-spanned by the remaining $d+2-k$ points, at an even number of points.

{\it Hint.} Perhaps one can take points on the moment curve.

(d'-k) No $d+1$ of $d+3$ given points in $\R^d$ lie in the same $(d-1)$-hyperplane.
Then the number of those $k$-simplices formed by $k+1$ of them that intersect at exactly one point the surface of the $(d+1-k)$-simplex formed by the remaining $d+2-k$ points, is even.

{\it Hint.} Induction on $k$; analogously to the alternative proof of Proposition \ref{r:lint3}.3'.
\end{conjecture}

\section{Realizability of products and the Menger conjecture}\label{s:prod}

\subsection{The Menger conjecture}\label{s:menger}

The {\it (Cartesian) product} $F\times F'$ of two figures $F,F'$ in $\R^3$ is the set of all points $(x,y,z,x',y',z')\in\R^6$ such that $(x,y,z)\in F$ and $(x',y',z')\in F'$.

\begin{figure}[h]
\centerline{\includegraphics[scale=1.05]{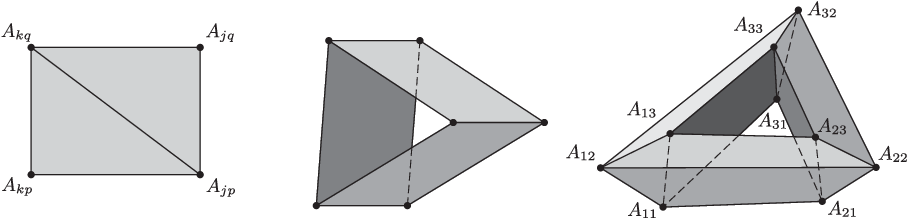}}
\caption{Realizations of the products: $K_2\times K_2$ (left), $K_2\times K_3$ (middle), $K_3\times K_3$ (right)}
\label{f-mnre}
\end{figure}


\begin{figure}[h]
\includegraphics[width=8cm]{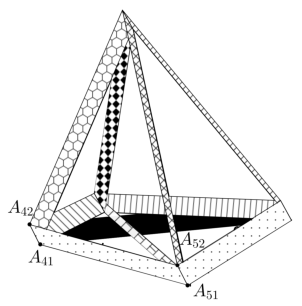}
\caption{Realization of the product $K_5\times K_2$}
\label{f:k5i}
\end{figure}

Examples of realization of products are given in Figures \ref{f-mnre} and \ref{f:k5i}.
For definition of realization see e.g. \cite[\S3.2]{Sk18}, \cite[\S6]{Sk}.
Karl Menger conjectured in 1929 that the square of a nonplanar graph is not realizable in $\R^4$ \cite{Me29} (cf. \aronly{Theorem \ref{0-ne4pr}}\jonly{\cite[Theorem 4.3]{Sk14}}).
This was proved only in 1978 by Brian Ummel \cite{Um78} using advanced algebraic topology.
A simple proof was obtained in 2003 by Mikhail Skopenkov \cite{Sk03} using \emph{lowering of dimension}, see exposition \aronly{below.}\jonly{in \cite[\S4]{Sk14}.}
His argument proves the generalized Menger conjecture (`the $k$-th power of a nonplanar graph is not realizable in $\R^{2k}$'), and even gives a short formula for the minimal number $d$ such that given product of several graphs is realizable in $\R^d$ \cite{Sk03}.

A combinatorial version of the product is \emph{product of two graphs} (not necessarily planar).
This product can be considered (although not canonically) as a hypergraph.
See Remark \ref{r:hyper}.a; cf. \cite[\S2, \S5]{ADN+}.

Propositions \ref{p:k4k4}.ab below imply that neither $K_4\times K_4$ nor $K_5\times K_3$ are (linearly) realizable in $\R^3$.
Proof of Proposition \ref{p:k4k4}.a shows that $K_{3,1}\times K_{3,1}$ is not realizable in $\R^3$ (but this product is realizable in $\R^4$).

\smallskip
{\bf Remark.}
Proofs of the Menger conjecture
using the van Kampen number or the Borsuk-Ulam theorem (see Remark \ref{s:compar}) are unknown.
The proof of the Menger conjecture in \cite{Um78} works for the topological version but is complicated.
The simpler proof in \cite{Sk03} uses for the topological version the non-trivial {\it Bryant approximation theorem}.
A simpler proof of the topological version can be obtained by inventing a quantitative PL version of the Menger conjecture (i.e., by improving the PL version of Theorem \ref{0-ne4pr} analogously to \S\ref{s:quant}, see
Remark \ref{c:menquan}).


\subsection{Realizability of products}\label{0statpr}


Let us formalize the idea of $K_m\times K_n$ drawn in 3- or 4-space.
Suppose that $A_{jp}$, where $j\in[m]$ and $p\in[n]$, are $mn$ points  in 3- or 4-space.
For numbers $j,k\in[m]$, $j<k$, and $p,q\in[n]$, $p<q$, denote
$$jk\times pq:=\{A_{jp}A_{kq}A_{jq},A_{jp}A_{kq}A_{kp}\}.$$
This is a pair of triangles having a common side (Figure \ref{f-mnre}, left).
Their union could be, but need not be, a plane quadrilateral.
An \textbf{$(m,n)$-product} is a collection of triangles from
$$jk\times pq,\quad\text{where}\quad 1\le j<k\le m,\quad 1\le p<q\le n.$$
(There are $mn(m-1)(n-1)/2$ such triangles.)
The \textbf{body} of $(m,n)$-product is the union of its triangles.
The body is a polyhedral and possibly self-intersecting

$\bullet$ square, if $m=n=2$ (Figure \ref{f-mnre},
left);

$\bullet$ lateral surface of a cylinder, if $m=3$ and $n=2$ (Figure \ref{f-mnre}, middle);

$\bullet$ torus, if $m=n=3$ (Figure \ref{f-mnre}, right).


\begin{proposition}\label{p:k4k4}
(a) Any $(4,4)$-product; \quad (b) Any $(3,5)$-product

in 3-space has a triangle and a side of a triangle which have disjoint vertices but intersect.
\end{proposition}


These folklore results are reduced to PL versions of Proposition \ref{0-k33}.ab in \S\ref{s:mnonre}.

\begin{theorem}[Square; \cite{Um78, Sk03}]\label{0-ne4pr}
Any $(5,5)$-product in 4-space has two triangles which have disjoint vertices but intersect.
\end{theorem}

The Square Theorem \ref{0-ne4pr} is reduced to
Theorem \ref{rlt-sachs} in \S\ref{s:mnre4}.

\begin{example}\label{stat-epr}
The analogues of Proposition \ref{p:k4k4} and Theorem~\ref{0-ne4pr} are false for

(a) $(2,n)$-products in 3-space for every $n$ (for $n\le 4$ this is obvious; for $n=5$ see Figure~\ref{f:k5i}:
the vertices of the parallelograms are the required $10$ points; for $n\ge6$ the construction is analogous, see \S\ref{s:mnre}; cf. \cite[Theorem 1.5]{RSS95'});

(b) $(3,n)$-products in 3-space for every $n\le4$ (for $n\le3$ this is obvious, see Figure~\ref{f-mnre}, right; for $n=4$ the construction is analogous, see \S\ref{s:mnre});

(c) $(4,n)$-products in 4-space for every $n$ (see \S\ref{s:mnre}).
\end{example}

\begin{remark}\label{c:menquan}
Denote
$$\t{K_5^2}\ :=\ \left\{\{(X,Y),(X',Y')\}\ :
\ X,Y,X',Y'\in{[5]\choose2},\text{ either } X\cap X'=\emptyset \text{ or }Y\cap Y'=\emptyset\right\}.$$
It would be interesting to find a subset $M\subset\t{K_5^2}$ such that
for any \emph{PL map} $f:K_5\times K_5\to\R^4$ there is an odd number of pairs $\{(X,Y),(X',Y')\}\in M$ for which $|f(X\times Y)\cap f(X'\times Y')|$ is odd, where by $X,Y,X',Y'$ we understand edges of $K_5$.
\footnote{For $M=\t{K_5^2}$ the number of such pairs has apparently the same parity as the number of linked unordered pairs of cycles in $K_{4,4}$ embedded in $\R^3$, which number is even.
So either $M$ should be different, or one should use integers (or residues modulo 4) instead of residues modulo 2.}

This is related to the following algebraic Menger problem \cite[Conjecture 2]{Pa20}:
{\it Complexes $K,L$ have non-trivial van Kampen obstructions to embeddability in $\R^m$ and in $\R^n$, respectively (see definition e.g. in \cite[\S1.5]{Sk18}).
Does the cartesian product $K\times L$ of $K$ and $L$ has non-trivial van Kampen obstruction to embeddability in $\R^{m+n}$?}
\end{remark}

\subsection{Realization of products in 3- and 4-space}\label{s:mnre}

\begin{proof}[Sketch of a proof of Example \ref{stat-epr}.a]
Let $A_{11},\ldots A_{1n}$ be points in $\R^3$ of which no four lie in one plane.
Take a vector $v$ not parallel to any plane passing through some three of these points.
For every $p\in[n]$ take a point $A_{2p}$ such that $\vec{A_{1p}A_{2p}}=v$.
If $v$ is small enough, then the points $A_{jp}$, $j\in\{1,2\}$, $p\in[n]$, are as required: there are no triangle and a side of a triangle with vertices at these points, which have disjoint vertices but intersect.

Indeed, $12\times pq$ is a parallelogram for every $p\ne q$.
Since no 4 of the points $A_{11},\ldots A_{1n}$ lie in one plane, and by the choice of $v$, for any distinct $p,q,r,s$ the segments $A_{1p}A_{1q}$ and $A_{1r}A_{1s}$ are disjoint.
Since $v$ is small enough, the same holds for 1 replaced by 2.
Then any two (convex hulls of) parallelograms $12\times pq$ and $12\times rs$ that have no common side are disjoint.
Now one can check that the points $A_{jp}$ are as required.
\end{proof}

\begin{proof}[Sketch of a proof of Example \ref{stat-epr}.b]
Let
$$A_{11}=(1,0,1),\quad A_{12}=(-1,0,1),\quad A_{13}=(0,0,2),\quad A_{14}=(0,0,3).$$
Let $f:\R^3\to \R^3$ be the rotation through $\frac{2\pi}{3}$ w.r.t. $x$-axis.
Let $A_{2p}=f(A_{1p})$ and $A_{3p}=f(f(A_{1p}))$ for every $p\in[4]$ (Figure \ref{f-mnre}, right).
Then the points $A_{jp}$, $j\in[3]$, $p\in[4]$, are as required.

Indeed, $jk\times pq$ is a parallelogram for every $j\ne k$, $p\ne q$.
Since every two segments joining points $A_{1p}$ either are disjoint or intersect at a common vertex,
any two of such parallelograms that have no common side are disjoint.
Now one can check that the points $A_{jp}$ are as required.
\end{proof}

\begin{figure}[h]
\includegraphics[width=6cm]{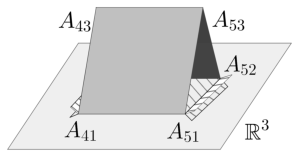}
\qquad\includegraphics[width=6cm]{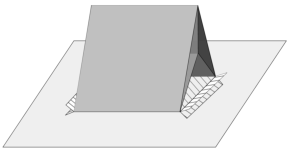}
\caption{Left: to realization in $\R^4$ of the product $K_3\times K_5$.
\newline
Right: to realization in $\R^4$ of the product $K_4\times K_5$.}
\label{f-k5xkn}
\end{figure}

\begin{proof}[Sketch of a proof of the weaker version of Example \ref{stat-epr}.c: $(3,5)$-product in 4-space]
Take a 3-dimensional hyperplane in $\R^4$ (shown in Figure \ref{f-k5xkn}, left, as a plane in 3-space).
In this hyperplane take 10 vertices $A_{jp}$, where $j\in[5]$, $p\in\{1,2\}$, shown in Figure~\ref{f:k5i}.
Take a vector $v$ not parallel to the hyperplane.
Set $A_{j3}:=A_{j1}+v$.
(In Figure \ref{f-k5xkn}, left, we see the lateral surface of the prismoid  $A_{41}A_{42}A_{43}A_{51}A_{52}A_{53}$.)
Then the points $A_{jp}$, $j\in[5]$, $p\in[3]$, are as required: there are no two triangles with vertices at these points, which have disjoint vertices but intersect.
\end{proof}

\begin{proof}[Sketch of a proof of Example \ref{stat-epr}.c]
Take points $A_{jp}\in \R^3\subset \R^4$, $j\in\{1,2\}$, $p\in[n]$ from the proof of Example \ref{stat-epr}.a.
Then $\overrightarrow{A_{1p}A_{1q}}=\overrightarrow{A_{2p}A_{2q}}$ for every $p\ne q$.
Take non-collinear vectors $v_3,v_4\in\R^4$ not parallel to the hyperplane $\R^3\subset \R^4$.
Denote $A_{jp}:=A_{1p}+v_j$, $j\in\{3,4\}$.
We could have taken $v_3,v_4$ so that $A_{14}$ is an interior point of the triangle $A_{11}A_{12}A_{13}$.
See Figure \ref{f-k5xkn}, right.
Then the points $A_{jp}$, $j\in[4]$, $p\in[n]$, are as required.
\end{proof}


\subsection{Non-realizability of products in 3-space}\label{s:mnonre}

\begin{figure}[h]\centering
\includegraphics{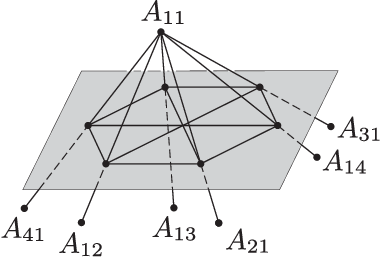}\qquad\begin{tikzpicture}[line cap=round,line join=round,>=triangle 45,x=0.32cm,y=0.32cm]
\clip(-4.91,-7.4) rectangle (15.45,4.56);
\fill[fill=black,fill opacity=0.1] (-3,-5) -- (0,1) -- (13,1) -- (10,-5) -- cycle;
\draw (-0.89,-3.76)-- (2.51,-3.71);
\draw (8.81,-3.6)-- (2.36,-3.71);
\draw (4,3)-- (0.89,-0.1);
\draw (4,3)-- (2.85,-1.32);
\draw (4,3)-- (-0.89,-3.76);
\draw (4,3)-- (10.95,0.08);
\draw (4,3)-- (8.81,-3.6);
\draw (4,3)-- (6.35,-2.52);
\draw (-1.79,-5)-- (-2.49,-5.96);
\draw [dash pattern=on 2pt off 2pt] (-0.89,-3.76)-- (-1.79,-5);
\draw [dash pattern=on 2pt off 2pt] (0.89,-0.1)-- (-2,-3);
\draw (-2,-3)-- (-3,-4);
\draw (1.87,-5)-- (1.39,-6.8);
\draw (9.83,-5)-- (10.48,-5.9);
\draw [dash pattern=on 2pt off 2pt] (1.87,-5)-- (2.85,-1.32);
\draw [dash pattern=on 2pt off 2pt] (8.81,-3.6)-- (9.83,-5);
\draw (-3,-5)-- (0,1);
\draw (13,1)-- (10,-5);
\draw (10,-5)-- (-3,-5);
\draw (0,1)-- (1.66,1);
\draw (3,1)-- (3.34,1);
\draw (4,1)-- (5,1);
\draw (6,1)-- (8,1);
\draw (9.32,1)-- (13,1);
\draw (2.85,-1.32)-- (-0.89,-3.76);
\draw (8.81,-3.6)-- (10.95,0.08);
\draw [dash pattern=on 2pt off 2pt] (6.35,-2.52)-- (7.4,-5);
\draw (7.4,-5)-- (8.15,-6.76);
\draw [dash pattern=on 2pt off 2pt] (10.95,0.08)-- (12.24,-0.51);
\draw (12.24,-0.51)-- (13.25,-0.96);
\draw (2.85,-1.32)-- (5.44,-0.88);
\draw (5.91,-0.8)-- (6.42,-0.71);
\draw (7.07,-0.6)-- (10.95,0.08);
\fill [color=black] (0.89,-0.1) circle (1.5pt);
\fill [color=black] (2.85,-1.32) circle (1.5pt);
\draw [color=black] (10.95,0.08) ++(-1.5pt,0 pt) -- ++(1.5pt,1.5pt)--++(1.5pt,-1.5pt)--++(-1.5pt,-1.5pt)--++(-1.5pt,1.5pt);
\fill [color=black] (8.81,-3.6) circle (1.5pt);
\fill [color=black] (6.35,-2.52) circle (1.5pt);
\draw [color=black] (-0.89,-3.76) ++(-1.5pt,0 pt) -- ++(1.5pt,1.5pt)--++(1.5pt,-1.5pt)--++(-1.5pt,-1.5pt)--++(-1.5pt,1.5pt);
\fill [color=black] (4,3) circle (1.5pt);
\draw[color=black] (3.07,3.76) node {$A_{11}$};
\fill [color=black] (-2.49,-5.96) circle (1.5pt);
\draw[color=black] (-1.13,-6.4) node {$A_{12}$};
\fill [color=black] (13.25,-0.96) circle (1.5pt);
\draw[color=black] (13.39,-1.94) node {$A_{13}$};
\fill [color=black] (-3,-4) circle (1.5pt);
\draw[color=black] (-3.21,-2.64) node {$A_{41}$};
\fill [color=black] (1.39,-6.8) circle (1.5pt);
\draw[color=black] (2.85,-6.64) node {$A_{21}$};
\fill [color=black] (10.48,-5.9) circle (1.5pt);
\draw[color=black] (11.53,-4.96) node {$A_{31}$};
\fill [color=black] (8.15,-6.76) circle (1.5pt);
\draw[color=black] (7.13,-6.56) node {$A_{51}$};
\end{tikzpicture}
\caption{To the proofs of Proposition \ref{p:k4k4}.a (left) and \ref{p:k4k4}.b (right)}
\label{f:3space}
\end{figure}

\begin{proof}[Proof of Proposition \ref{p:k4k4}.a]
(The proof is analogous to Proposition~\ref{0-ne3}.)
There is a plane $\alpha$ such that one of the given points lies on one side of $\alpha$, and the other 15 points lie on the other side of $\alpha$ (Figure~\ref{f:3space}, left).
The intersection point of $\alpha$

$\bullet$ with the segment $A_{jp}A_{kp}$ is colored in blue for every $k\in[4]-\{j\}$;

$\bullet$ with the segment $A_{jp}A_{jq}$ is colored in red for every $q\in[4]-\{p\}$.

The intersection of a product $jk\times pq$ with $\alpha$ is called an {\it arc}.
Then arcs have ends of different color.
(The intersection of $\alpha$ with the body of the $(4,4)$-product is a PL drawing, possibly with self-intersections, of $K_{3,3}$ in $\alpha$, i.e., the image of a {\it PL map} $K_{3,3}\to\alpha$.)

Then by the PL analogue of Proposition \ref{0-k33}.a\footnote{This analog is proved analogously, see  \cite[Remark 1.4.4b]{Sk18}.} there are intersecting arcs without common edges.
Then there are $k,k'\in[4]-\{j\}$ and $q,q'\in[4]-\{p\}$ such that $k\ne k'$, $q\ne q'$, and some two triangles, one from $jk\times pq$ and the other from $jk'\times pq'$, have a common point distinct from the common vertex  $A_{jp}$ of the triangles.
Hence one of these triangles intersects a side of the other not passing through $A_{jp}$, hence not having any common vertices with the first triangle.
\end{proof}



Given 9 points $A_{jp}$, $j\in\{u,v,w\}$, $p\in\{u',v',w'\}$, in 3- or in 4-space denote by $uvw\times u'v'w'$ the corresponding $(3,3)$-product (Figure \ref{f-mnre}, right; as opposed to the figure, the $(3,3)$-product can have self-intersections).


\begin{proof}[Proof of Proposition \ref{p:k4k4}.b]
We may assume that no 4 of given $15$ points $A_{jp}$, $j\in[3]$, $p\in[5]$, lie in one hyperplane
(analogously to Remark~\ref{r:stronger}.b).
There is a plane $\alpha$ such that one of the given points, say $A_{j,p}$, lies on one side of $\alpha$, and the other 15 points lie on the other side of $\alpha$ (Figure~\ref{f:3space}, right).
The intersection point of $\alpha$

$\bullet$ with the segment $A_{jp}A_{kp}$ is colored in blue for every $k\in[3]-\{j\}$;

$\bullet$ with the segment $A_{jp}A_{jq}$ is colored in red for every $q\in[5]-\{p\}$.

The intersection of a product $jk\times pq$ with $\alpha$ is called an {\it arc}.
Then arcs have ends of different color.
(The intersection of $\alpha$ with the union of the triangles of the $(5,3)$-product is
a PL drawing, possibly with self-intersections, of $K_{2,4}$ in $\alpha$.)

Analogously to the last paragraph of the proof of Proposition \ref{p:k4k4}.a either

(1) the $(5,3)$-product has a triangle and a side of a triangle which have disjoint vertices but intersect, or

(2) any two arcs intersect only at their common vertex (if there is one).

In the second case denote the blue points by $K,L$.
By the PL analogue of Proposition \ref{0-k33}.b there are two red points $Q,R$ such that the remaining two red points $S,T$ lie on different sides w.r.t. the closed polygonal $\gamma$ formed by arcs $QK, KR, RL, LQ$.
Take $k,l\in[3]-\{j\}$ and $q,r,s,t\in[5]-\{p\}$ such that the points $K,L$, and $Q,R,S,T$ belong to the segments joining $A_{jp}$ to $A_{kp},A_{lp}$, and to $A_{jq},A_{jr},A_{js},A_{jt}$, respectively.
Then $\alpha$ intersects

$\bullet$ the outline of the triangle $j\times pst:=A_{jp}A_{js}A_{jt}$ by $S$ and $T$ (note that the triangle is not contained in the $(3,5)$-product);

$\bullet$ the body of the $(3,3)$-subproduct $jkl\times pqr$ (contained in the $(3,5)$-product) by $\gamma$.

If $\gamma$ has self-intersections, then we obtain the property (1).
If not, then we obtain the property (1) by Lemma \ref{0-even3} below because $S,T$ lie on different sides w.r.t. $\gamma$.
\end{proof}

\begin{lemma}[lowering of dimension]\label{0-even3}
In 3-space the outline $\partial\Delta$ of a triangle $\Delta$ and a $(3,3)$-product $\tau$ have a unique common vertex $O$.
No 4 vertices of $\Delta$ and $\tau$ lie in one plane.
A plane $\alpha$ splits $O$ from the base of $\Delta$, and from the remaining 8 points of $\tau$.
The plane $\alpha$ intersects

$\bullet$ the body $|\tau|$ by a closed polygonal line without self-intersections;

$\bullet$ $\partial\Delta$ by two points $S,T$ lying on different sides w.r.t. the polygonal line.

Then $\tau$ and $\partial\Delta$ have a triangle and a side of a triangle which have disjoint vertices but intersect.
\end{lemma}

\begin{proof}
Denote by $\left<\Delta\right>$ the plane of $\Delta$.
Then $\alpha\cap\left<\Delta\right>$ is a line.
The intersection $\alpha\cap\Delta$ is the segment $ST$.
The points $S,T$ lie in $\alpha$ on different sides w.r.t. the closed polygonal line $\alpha\cap|\tau|$ which does not have self-intersections, and no 3 points among $S,T$ and the vertices of $\alpha\cap|\tau|$ lie in one line.
Hence the segment $\alpha\cap\Delta$ intersects $|\tau|\cap\left<\Delta\right>$ at an odd number of points.
Analogously, the line $\alpha\cap\left<\Delta\right>$ intersects $|\tau|\cap\left<\Delta\right>$ at an even  number of points.
Since also no 4 vertices of $\Delta$ and $\tau$ lie in one plane, it follows that $\tau\cap\left<\Delta\right>:=\{\Gamma\cap\left<\Delta\right>\ :\ \Gamma\in\tau\}$ is
a \emph{1-cycle}, i.e., is a set of segments in the plane $\left<\Delta\right>$ such that every point of
$\left<\Delta\right>$ is the endpoint of an even number (possibly, zero) of the segments.
By the assumption on $\alpha$ we may choose the segments so that the line $\alpha\cap\left<\Delta\right>$ splits 
$O$ from $S,T$, and from all the vertices of the segments distinct from $O$ (cf. Figure \ref{f:plane}, left).
Hence by an analogue of Lemma \ref{0-eveng} for a triangle and a 1-cycle, some two segments from $\partial\Delta$ and from $\tau\cap\left<\Delta\right>$ intersect but do not have common vertices.
At most one of these segments contains $O$.
Hence the obtained segment of the 1-cycle $\tau\cap\left<\Delta\right>$ is the intersection with  $\left<\Delta\right>$ of a triangle from $\tau$, which intersects a side of $\Delta$, but does not have common vertices with the side.
\end{proof}


\subsection{Non-realizability of products in 4-space}\label{s:mnre4}

A \emph{Seifert chain} (or a coboundary) of a closed polygonal line $a$ in 3-space is a finite collection $S$ of triangles (non-degenerate to a segment or a point) in 3-space such that

$\bullet$ every edge of $a$ is the side of exactly one triangle from $S$;

$\bullet$ every segment that is not an edge of $a$ is the side of an even number (possibly, zero) of triangles from $S$.

Two disjoint closed polygonal lines $a$ and $a'$ in 3-space \emph{linked modulo 2} if for any Seifert chains $S$ of $a$ and $S'$ of $a'$ such that the outline of any triangle of $S$ is disjoint from the outline of any triangle of $S'$, the  number of linked modulo 2 pairs $(\Delta,\Delta')$ of triangles $\Delta$ of $S$ and $\Delta'$ of $S'$ is odd.
The equivalence to other definitions of being linked modulo 2 (in particular, to the definition before Theorem \ref{rlt-sachs}) is proved in \cite[Proposition 4.8.3]{Sk}.

\begin{proof}[Proof of  the Square Theorem \ref{0-ne4pr}]
We may assume that no 5 of the given 25 points $A_{jp}$, $j,p\in[5]$, lie in one 3-dimensional hyperplane
(analogously to Remark \ref{r:stronger}.b).
There is a 3-dimensional hyperplane $\alpha$ such that one of the given points, say $A_{j,p}$, lies on one side of $\alpha$, and the other 24 points lie on the other side of $\alpha$ (Figure~\ref{f:4spacepr}).
The intersection point of $\alpha$

$\bullet$ with the segment $A_{jp}A_{kp}$ is colored in blue for every $k\in[5]-\{j\}$;

$\bullet$ with the segment $A_{jp}A_{jq}$ is colored in red for every $q\in[5]-\{p\}$.

The intersection of a product $jk\times pq$ with $\alpha$ is called an {\it arc}.
Then arcs have ends of different color.

\begin{figure}[h]\centering
\begin{tikzpicture}[line cap=round,line join=round,>=triangle 45,x=0.32cm,y=0.32cm]
\clip(-4.93,-7.54) rectangle (15.43,4.42);
\fill[fill=black,fill opacity=0.1] (-3,-5) -- (0,1) -- (13,1) -- (10,-5) -- cycle;
\draw (-1,-4)-- (7,-2);
\draw (3.03,-1.9)-- (4.53,-2.43);
\draw (9,-4)-- (5.03,-2.61);
\draw (7,-2)-- (5.54,-1.51);
\draw (4,3)-- (1,0);
\draw (4,3)-- (-1,-4);
\draw (4,3)-- (3.03,-1.9);
\draw (4,3)-- (7,-2);
\draw (4,3)-- (9,-4);
\draw (4,3)-- (11,0);
\draw (2.39,-5)-- (2.03,-6.74);
\draw [dash pattern=on 2pt off 2pt] (3.03,-1.9)-- (2.39,-5);
\draw [dash pattern=on 2pt off 2pt] (7,-2)-- (8.65,-5);
\draw (8.65,-5)-- (9.47,-6.48);
\draw [dash pattern=on 2pt off 2pt] (1,0)-- (-2,-3);
\draw (-2,-3)-- (-3,-4);
\draw (-1.71,-5)-- (-2.4,-5.96);
\draw (9.71,-5)-- (10.61,-6.25);
\draw (12.24,-0.53)-- (13.38,-1.02);
\draw [dash pattern=on 2pt off 2pt] (-1.71,-5)-- (-1,-4);
\draw [dash pattern=on 2pt off 2pt] (9,-4)-- (9.71,-5);
\draw [dash pattern=on 2pt off 2pt] (11,0)-- (12.24,-0.53);
\draw (-3,-5)-- (0,1);
\draw (13,1)-- (10,-5);
\draw (10,-5)-- (-3,-5);
\draw (0,1)-- (1.33,1);
\draw (3,1)-- (3.34,1);
\draw (4,1)-- (5,1);
\draw (6,1)-- (8,1);
\draw (9.32,1)-- (13,1);
\draw (3.03,-1.9)-- (6.27,-1.13);
\draw (7.98,-0.72)-- (11,0);
\draw (1,0)-- (1.5,-0.17);
\draw (1.89,-0.3)-- (2.85,-0.62);
\draw (3.57,-0.86)-- (5.02,-1.34);
\draw (1,0)-- (-0.45,-0.72);
\draw (-0.45,-0.72)-- (-1,-4);
\draw (11,0)-- (10.01,-3.36);
\draw (10.01,-3.36)-- (9,-4);
\draw (4,3)-- (-0.45,-0.72);
\draw (4,3)-- (10.01,-3.36);
\draw [dash pattern=on 2pt off 2pt] (-0.45,-0.72)-- (-1.16,-1.32);
\draw [dash pattern=on 2pt off 2pt] (10.01,-3.36)-- (10.54,-3.92);
\draw (-1.16,-1.32)-- (-2.3,-2.27);
\draw (10.54,-3.92)-- (11.51,-4.96);
\draw [color=black] (1,0) ++(-1.5pt,0 pt) -- ++(1.5pt,1.5pt)--++(1.5pt,-1.5pt)--++(-1.5pt,-1.5pt)--++(-1.5pt,1.5pt);
\draw [color=black] (-1,-4) ++(-1.5pt,0 pt) -- ++(1.5pt,1.5pt)--++(1.5pt,-1.5pt)--++(-1.5pt,-1.5pt)--++(-1.5pt,1.5pt);
\fill [color=black] (7,-2) circle (1.5pt);
\draw [color=black] (9,-4) ++(-1.5pt,0 pt) -- ++(1.5pt,1.5pt)--++(1.5pt,-1.5pt)--++(-1.5pt,-1.5pt)--++(-1.5pt,1.5pt);
\draw [color=black] (11,0) ++(-1.5pt,0 pt) -- ++(1.5pt,1.5pt)--++(1.5pt,-1.5pt)--++(-1.5pt,-1.5pt)--++(-1.5pt,1.5pt);
\fill [color=black] (3.03,-1.9) circle (1.5pt);
\fill [color=black] (4,3) circle (1.5pt);
\draw[color=black] (3.05,3.72) node {$A_{11}$};
\fill [color=black] (2.03,-6.74) circle (1.5pt);
\draw[color=black] (3.45,-6.26) node {$A_{14}$};
\fill [color=black] (9.47,-6.48) circle (1.5pt);
\draw[color=black] (8.45,-6.34) node {$A_{13}$};
\fill [color=black] (-3,-4) circle (1.5pt);
\draw[color=black] (-3.91,-3.16) node {$A_{21}$};
\fill [color=black] (-2.4,-5.96) circle (1.5pt);
\draw[color=black] (-0.93,-5.8) node {$A_{31}$};
\fill [color=black] (10.61,-6.25) circle (1.5pt);
\draw[color=black] (12.35,-6.36) node {$A_{41}$};
\fill [color=black] (13.38,-1.02) circle (1.5pt);
\draw[color=black] (13.21,-2.14) node {$A_{51}$};
\fill [color=black] (-0.45,-0.72) circle (1.5pt);
\fill [color=black] (10.01,-3.36) circle (1.5pt);
\fill [color=black] (-2.3,-2.27) circle (1.5pt);
\draw[color=black] (-2.93,-0.94) node {$A_{12}$};
\fill [color=black] (11.51,-4.96) circle (1.5pt);
\draw[color=black] (13.41,-4.8) node {$A_{15}$};
\end{tikzpicture}
\caption{To the proof of the Square Theorem~\ref{0-ne4pr}}
\label{f:4spacepr}
\end{figure}


Analogously to the last paragraph of the proof of Proposition \ref{p:k4k4}.a either

(1) the $(5,5)$-product has two triangles which have disjoint vertices but intersect, or

(2) any two arcs intersect only at their common vertex (if there is one).

In the second case the intersection of $\alpha$ with the body of the $(5,5)$-product is a PL drawing without self-intersections of $K_{4,4}$ in $\alpha$.
Use the following PL analogue of Theorem~\ref{rlt-sachs}: in any PL drawing without self-intersections of $K_{4,4}$ in $\R^3$ there are two cycles of length 4 which are linked modulo 2 (see proof in \cite{Sa81, Z13}).
We obtain two linked modulo 2 closed polygonal lines in $\alpha$, each consisting of four arcs.
Take $\{a,b,a',b'\}=[5]-\{j\}$ and $\{c,d,c',d'\}=[5]-\{p\}$ such that the arcs

$\bullet$ of the first polygonal line belong to the products $ja\times pc$, $jb\times pc$, $ja\times pd$, $jb\times pd$,

$\bullet$ of the second polygonal line belong to the products $ja'\times pc'$, $jb'\times pc'$, $ja'\times pd'$, $jb'\times pd'$.

Then the polygonal lines are the intersections with the hyperplane of the bodies of the $(3,3)$-products   $jab\times pcd$ and $ja'b'\times pc'd'$.
So the required statement is implied by the following Lemma~\ref{0-even4g}.
\end{proof}

\begin{lemma}[lowering of dimension]\label{0-even4g}
Two $(3,3)$-products in 4-space have a unique common vertex $O$.
No 5 of their $17$ vertices lie in one 3-dimensional hyperplane.
A 3-dimensional hyperplane $\alpha$ splits
$O$ from the remaining 16 points of the $(3,3)$-products.
The hyperplane $\alpha$ intersects the bodies of the $(3,3)$-products by a pair of disjoint closed polygonal lines linked modulo 2 in $\alpha$.
Then the $(3,3)$-products have two triangles which have disjoint vertices but intersect.
\end{lemma}


\begin{proof}
Denote by $S$ (by $S'$) the set of all triangles from the first (the second) $(3,3)$-product, which do not contain $O$.
Since no $5$ of the $17$ vertices of the $(3,3)$-products lie in one 3-dimensional hyperplane, no $4$ of their $16$ projections to $\alpha$ with the center $O$ lie in one plane.
Then the outlines of the triangles $\alpha\cap O\Delta$ and $\alpha\cap O\Delta'$ are disjoint
for any $\Delta\in S$ and $\Delta'\in S'$.
Denote by $\gamma$ and $\gamma'$ the given disjoint closed polygonal lines.
Since $\gamma$ and $\gamma'$ are linked modulo 2 in $\alpha$, the number of linked modulo 2 pairs $(\Delta,\Delta')$ of such triangles is odd.
By Lemma \ref{relation}.a such triangles are linked modulo 2 if and only if the surfaces of the tetrahedra   $O\Delta$ and $O\Delta'$ intersect at an even number of points (including $O$).
Take any side $MN$ of a triangle from $S$, and any side $M'N'$ of a triangle from $S'$.
If $MN$ is not contained in $\gamma$, and $M'N'$ is not contained in $\gamma'$, the intersection $OMN\cap OM'N'$ appears exactly in two intersections of lateral surfaces of tetrahedra $O\Delta$ and  $O\Delta'$.
Then the numbers of pairs of intersecting triangles having one of the following types is odd:

$\bullet$ pairs $(\Delta,\Delta')$ for $\Delta\in S$ and $\Delta'\in S'$;

$\bullet$ pairs $(OMN,\Delta')$ for a side $MN$ of $\gamma$, and a triangle $\Delta'\in S'$.

$\bullet$ pairs $(\Delta,OM'N')$ for a side $M'N'$ of $\gamma'$, and a triangle $\Delta\in S$.

Now the lemma follows because in any of these pairs triangles have disjoint vertices, and are contained in triangles of given $(3,3)$-products.
\end{proof}

\comment

Denote by $\Delta_1,\ldots,\Delta_9$ ($\Delta_1',\ldots,\Delta_9'$) those triangles of $T$ (of $T'$) that do not contain $O$.
Let $OX=\conv\{\{O\}\cup X\}$ be the cone over $X$ with the center $O$.
Then $(T\cap T')-\{O\}=\emptyset$ consists of an even number of points.
Hence there is an even number of pairs $(j,j')\in[9]^2$ such that the surfaces of tetrahedra $O\Delta_j$ and $O\Delta_{j'}'$ intersect at an odd number of points.
By (a spherical analogue of) the argument in the proof of Theorem \ref{0-ne4}
the latter number has the same parity as the number of pairs $(j,j')\in[9]^2$ such that the triangles $\pi\cap O\Delta_j$ and $\pi\cap O\Delta_{j'}'$ are linked.
So the lemma follows by (a spherical analogue of) Proposition \ref{p:linqua}.b.

partly to algolink, algor



\subsection{Parity Lemmas}\label{s:par3}

!!!

For the proof of the Square Theorem \ref{0-ne4pr}  we need Lemma \ref{0-even4g} whose simpler analogues  were already used above (see Lemma \ref{0-eveng} and the argument on a triangle and a $(3,3)$-product in 3-space from the proof of \ref{p:k4k4}.b).
See a similar alternative proof in \cite{Zu} and more on parity lemmas in \cite[\S1.3]{Sk18}, \cite[\S4]{Sk}.

In order to illustrate the idea in simpler situations, we start with a 2-dimensional Parity Lemma \ref{0-even}, then proceed to a 3-dimensional Parity Lemma \ref{rlt-lk4}, and then present Proposition \ref{p:linqua} on linking in 3-space.
All of them are required for Lemma \ref{0-even4g}.

Further we generalize the following evident fact:
{\it if no 4 of the vertices of a closed polygonal line and of a tetrahedron in 3-space  lie in one plane, then the polygonal line and the surface of the tetrahedron intersect at an even number of points.}


Some points in 3-space {\bf are in general position}, if

$\bullet$ no 4 of them  lie in one plane, and

$\bullet$ for every pairwise disjoint pair, triple and triple among them the line passing through the pair, the plane passing through the first triple and  the plane passing through the second triple, have no common points.



A {\bf 2-(dimensional) cycle (modulo 2)} in 3-space is a set of triangles (non-degenerate to a segment or to a point) in 3-space such that every segment in 3-space is a side of an even number of triangles from the set.
An example of a 2-cycle is the collection of faces of a tetrahedron (or of another `triangulated surface').
{\it The vertices} of a 2-cycle are the vertices of its triangles.
{\it The body} of a 2-cycle is the union of its triangles.
Cf. \cite[\S3]{ADN+}.

An example of a 2-cycle is the surface of a tetrahedron (possibly, degenerate).
Also, (the body of) the $(3,3)$-product $T_{uvw}$  defined in \S\ref{s:mnonre} is the body of a 2-cycle.



\begin{lemma}[Parity]\label{rlt-lk4}
If the vertices of a closed polygonal line and a 2-cycle in 3-space are in general position, then the polygonal line intersects the body of the 2-cycle at an even number of points.
\end{lemma}

\footnote{It is here that we use a specific triangulation of $K_4\times K_4$.
Thus the point $A_{11}$ is not interchangeable with other $A_{jp}$.
So we have to consider a tetrahedron instead of a (hyper)plane as in Theorems \ref{rlt-cgslin}, \ref{0-ne4} and \ref{0-ne4pr}.
Analogous remark applies for the proof of \ref{p:k4k4}.b below. }

Lemma \ref{0-even3} is an analogue of Lemma \ref{0-eveng}.
Unlike those lemmas, Lemma \ref{0-even3} is stated in the `qualitative' not `quantitative' way because we do not prove the `quantitative' version of \ref{p:k4k4}.b, and because the `quantitative' way  requires general position assumption (see definition in \S\ref{s:par3}).


The following analogue of Lemma \ref{0-eveng} is implied by its obvious particular case when $\gamma$ is just one closed polygonal line.
(An accurate proof is analogous to accurate proof of Lemma \ref{0-eveng} in \S\ref{s:par3}.)

{\it Let $\Delta$ be a triangle and $\gamma$ the union of a finite number of closed polygonal lines
(not necessarily disjoint) in the plane.
Assume that $\Delta$ and $\gamma$ intersect at a unique point, which is their common vertex,
and which is contained in the interior of a small polygon.
Suppose that the outline of the polygon intersects $\Delta$
and $\gamma$ by a polygonal line $\delta$ and by a finite collection of points, respectively.
Then the number of intersection points of $\delta$ and $\gamma$ is even.}


Lemma \ref{0-even3} is reduced to this analogue by proving that the intersection of the $(3,3)$-product $T_{123}$ and the plane containing the triangle is a finite union of closed polygonal lines.



\begin{proof}[Proof of Lemma \ref{0-eveng}]
Denote the point by $O$ and the triangles by $OX'Y'$ and $OZ'T'$ so that $X,Y,Z,T$ are the intersection points of the line and  $OX',OY',OZ',OT'$ respectively.
Let $a:=\partial(OX'Y')$ and $b:=\partial(OZ'T')$.
Take points $X_1,Y_1$ in the line so close to $X,Y$ that $(\partial(X'XX_1)\cup\partial(Y'YY_1))\cap b=\emptyset$.
The lemma follows because
$$|XY\cap\{Z,T\}|=|XY\cap b|=|\partial(OXY)\cap b|-1 \underset2\equiv
|a\cap b|+|\partial(XYY'X')\cap b|-1\underset2\equiv |a\cap b|-1.$$
Here the last congruence holds because
$|\partial(XYY'X')\cap b|=|\partial(X_1Y_1Y'X')\cap b|\underset2\equiv0$ by the Parity Lemma \ref{0-even}
\end{proof}

\begin{proposition}\label{p:linqua} Let $ABCD$ and $A'B'C'D'$ be two closed quadrangular polygonal lines in 3-space no 4 of whose 8 vertices lie in one  plane.

(a) The polygonal lines are linked if and only if an odd number among the following pairs of triangles are linked pairs:
$$(ABC,A'B'C'),\quad (ABC,A'D'C'),\quad(ADC,A'B'C'),\quad(ADC,A'D'C').$$

(b) Assume that $\Delta_1,\ldots,\Delta_k$ are triangles in 3-space such that $\Delta_1,\ldots,\Delta_k,ABC,ADC$ is a 2-cycle and the union of their vertices is in general position.
(Such a collection of triangles is called a {\it coboundary} of $ABCD$.)
Assume that $\Delta_1',\ldots,\Delta_{k'}'$ is an analogous collection of triangles for $A'B'C'D'$.
The polygonal lines are linked if and only if an odd number among the $kk'$ pairs $(\Delta_j,\Delta_{j'}')$ of triangles are linked pairs.
\end{proposition}

\begin{proof} Part (a) is a particular case of (b) for $k=k'=2$, $\Delta_1=ABC$, $\Delta_2=ADC$, $\Delta_1'=A'B'C'$, $\Delta_2'=A'D'C'$.

Denote by $\partial\Delta$ the outline of a triangle or a quadrilateral $\Delta$.
Part (b) follows because
$$|ABCD\cap(A'B'C'\cup A'D'C')| \underset2\equiv \sum\limits_{j'=1}^{k'}|ABCD\cap\Delta_{j'}'| \underset2\equiv
\sum\limits_{j=1,\ j'=1}^{k,k'}|(\partial\Delta_j)\cap\Delta_{j'}'|.$$
Here the first congruence follows by the Parity Lemma \ref{rlt-lk4}.
 \end{proof}

\endcomment


{\it Books, surveys and expository papers in this list are marked by the stars.}

\end{document}


\newcommand{\mytitle}{Реализуемость гиперграфов и неотъемлемая зацепленность}
\title[\mytitle]{\mytitle}

\author{А. Б. Скопенков}

\thanks{Домашняя страница: \url{https://users.mccme.ru/skopenko/}.
\newline
Этот обзор основан на лекциях автора в Московском Физико-Техническом Инститите, Независимом Московском Университете, Институте Науки и Технологии (Австрия), и в различных летних школах и математических кружках.
Благодарю за перевод текста на русский язык С.~Дженжера, а также за полезные обсуждения О.~Виро, С.~Дженжера, А.~Зимина, Дж.~Зунга, И.~Изместьева, Р.~Карасева, А.~Руховича, М.~Скопенкова, М. Старкова, А.~Сосинского, С.~Табачникова, Г.~Челнокова,  А.~Шаповалова и анонимного рецензента от сборника <<Мат. Просвещение>>.}

\date{}
\maketitle

\begin{abstract}
    Мы приводим простые доказательства теоремы Конвея-Гордона-Закса о неотъемлемой зацепленности в трехмерном пространстве, а также теорем ван Кампена-Флореса и Уммеля о неотъемлемых пересечениях.
    Последние связаны с нереализуемостью некоторых гиперграфов в четырехмерном пространстве.
    Доказательства используют сведение к меньшим размерностям, которое позволяет показать связь между этими результатами.
    Используя элементарный язык, мы преподносим основные идеи без технических подробностей.
    Поэтому изложение доступно неспециалистам в данной области, включая студентов и школьников, знающих основы стереометрии и готовых изучить простые четырехмерные обобщения.
\end{abstract}

\tableofcontents

\bigskip
\newpage
\hfill{\it `It's too difficult.'}

\hfill{\it `Write simply.'}

\hfill{\it `That's hardest of all.'}

\hfill{\it I. Murdoch, The Message to the Planet.}

\section{Введение}\label{0intr}

\subsection{Невозможные конструкции, неотъемлемые пересечения и зацепленность}

\begin{figure}[h]
 \includegraphics[scale=0.7]{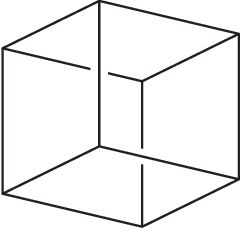} \qquad  \includegraphics[scale=0.6]{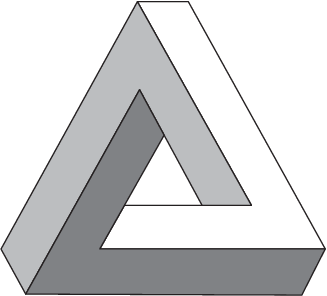} \qquad  \includegraphics[scale=0.7]{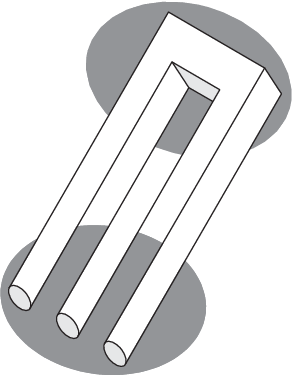} \qquad  \includegraphics[scale=0.9]{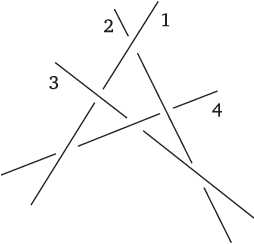}
\caption{Невозможный куб, треугольник Пенроуза, трезубец, невозможная проекция}
\label{f-imp}
\end{figure}

<<Невозможные конструкции>>, такие как невозможный куб, треугольник Пенроуза, невозможный трезубец и т.~д., хорошо известны, в основном благодаря рисункам Маурица Корнелиса Эшера, см. рис.~\ref{f-imp} и \cite{Io}, а также \cite{Br68, CKS+, GSS+}.
Эти изображения не допускают глобальной пространственной реализации из-за несоответствия локальных пространственных реализаций друг другу.
В геометрии, топологии и теории графов также известны базовые примеры <<невозможных конструкций>>
(в которых локальные части <<возможны>>).

 \begin{figure}[h]\centering
 \includegraphics{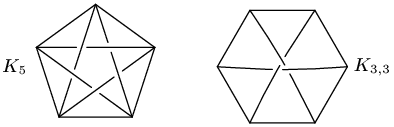}
 \caption{Непланарные графы $K_5$ и $K_{3,3}$}
 \label{k5}
 \end{figure}

Следующий пример <<невозможной конструкции>> или <<неотъемлемого пересечения>> уже связан с данной статьей. 
\emph{Ни для каких пяти точек на плоскости невозможно попарно соединить их путями так, чтобы пути пересекались только в своих начальных или конечных точках.}\footnote{Также невозможно расставить три домика и три колодца на плоскости и соединить путем каждый домик с каждым колодцем  так, чтобы пути пересекались только в своих начальных или конечных точках.
В терминах теории графов эти утверждения означают, что ни полный граф $K_5$ на пяти вершинах, ни полный двудольный граф $K_{3,3}$ не планарен, см. рис.~\ref{k5}.
Нижеприведенное утверждение \ref{0-ra2} является <<линейной>> версией непланарности графа $K_5$.}

\begin{proposition}[теорема о неотъемлемых пересечениях в размерности 2]\label{0-ra2} 
Для любых $5$ точек на плоскости найдутся два пересекающиеся отрезка с концами в этих точках, не имеющие общих вершин.
\end{proposition}


Далее понадобятся следующие соглашения.
Мы сокращаем <<трехмерное евклидово пространство $\R^3$>> до <<пространство>>, а
<<$d$-мерное евклидово пространство $\R^d$>> до <<$d$-мерное пространство>> или до <<$\R^d$>> (в основном $d=4$).
\textbf{Треугольником} называется часть плоскости, ограниченная замкнутой трехзвенной ломаной.

\begin{figure}[h]\centering
\includegraphics{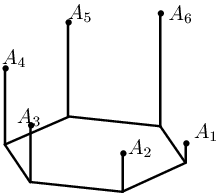} \qquad  \includegraphics{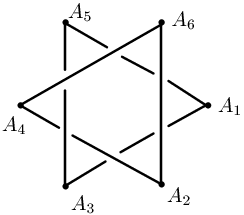} \qquad \includegraphics{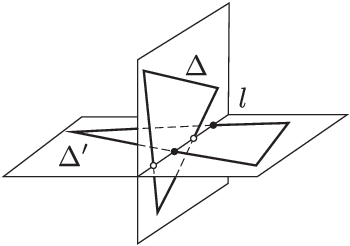}
\caption{Зацепленные треугольники}
\label{genpos}
\end{figure}

Возьмем два треугольника в пространстве, никакие $4$ из $6$ вершин которых не лежат в одной плоскости.
Эти треугольники называются \textbf{зацепленными}, если контур первого из них пересекает второй треугольник ровно в одной точке.
Например, треугольники $A_1A_3A_5$ и $A_2A_4A_6$, $\Delta$ и $\Delta'$ на рис.~\ref{genpos} зацеплены.\footnote{Расстояние от точки $A_j$ до плоскости проекции равно $j$.
Поэтому проекция на рис.~\ref{genpos}, посередине, реализуема, в отличие от рис.~\ref{f-imp} справа.
\newline
Свойство зацепленности симметрично (это не очевидно из определения, но имеет простое доказательство \cite[задача 4.1.1]{Sk}).
Другие вводные результаты о зацепленных треугольниках приведены в \cite[\S4.1 <<Зацепленность треугольников>>]{Sk}, \cite[\S1.1]{Sk24}.}

\begin{theorem}[о неотъемлемой зацепленности в размерности 3; Конвей--Гордон--Закс \cite{Sa81, CG83}]\label{rlt-cgslin}
Если никакие $4$ из данных $6$ точек в пространстве не лежат в одной плоскости, то найдутся два зацепленных треугольника с вершинами в данных $6$ точках.
\end{theorem}


\begin{theorem}[о неотъемлемых пересечениях в размерности 4; ван Кампен-Флорес \cite{vK32, Fl34}]\label{0-ne4}
Для любых $7$ точек в четырехмерном пространстве найдутся два пересекающиеся треугольника с вершинами в этих точках, не имеющие общих вершин.
\end{theorem}
 

 \begin{figure}[h]
 \includegraphics[width=6.5cm]{real-4.mps}
 \caption{Пять точек в $\R^3$ (реализация полного $3$-однородного гиперграфа на $5$ вершинах)}
 \label{f-full5}
 \end{figure}

Аналог теорем в размерности 2 и 4 верен для $6$ точек в пространстве (утверждение \ref{0-ne3}.b).
Аналоги теорем в размерности 2 и 4 неверны для $4$ точек на плоскости и для $6$ точек в четырехмерном пространстве, соответственно. 
Действительно (ср. рис.~\ref{f-full5}), возьмем 

$\bullet$ на плоскости три вершины и внутреннюю точку треугольника; 

$\bullet$ в четырехмерном пространстве пять вершин и внутреннюю точку четырехмерного симплекса.


\begin{remark}[о понижении размерности]\label{s:lowering}
Часто удобно свести результат на плоскости к результату на прямой, а трехмерный результат к результату на плоскости.
Аналогично четырехмерный результат можно сводить к трехмерному. 
См.  примеры в начале~\S\ref{s:interest}.
Благодаря такому <<понижению размерности>> читатель, незнакомый с четырехмерным пространством, не должен пугаться.
Интуиции о четырехмерном пространстве, выходящем за пределы замечания \ref{s:intu4}, не требуется для доказательства теоремы \ref{0-ne4} в размерности 4. 
\end{remark}

\subsection{Почему этот обзор может быть интересным?}\label{s:interest}

Мы описываем связь между неотъемлемой зацепленностью и неотъемлемыми пересечениями в последовательных размерностях. 
Она обобщается на высшие размерности (см. теорему~\ref{stat-il}; мы упоминаем размерности выше четырех только в этой теореме). 



Более конкретно, теорема в размерности 2 (утверждение~\ref{0-ra2}) сводится (в \S\ref{0pla}) к следующей теореме в размерности 1.  

\begin{proposition}[очевидное; теорема о неотъемлемой зацепленности в размерности 1]\label{p:liline}
Из любых $4$ точек на прямой можно покрасить две в красный, а две в синий цвет так, что цвета будут чередоваться: красный-синий-красный-синий или синий-красный-синий-красный (иными словами, <<красная пара зацеплена с синей>>).
\end{proposition}

Аналогично теорема~\ref{rlt-cgslin} в размерности 3 сводится (в \S\ref{0ramcgs}) к теореме в размерности 2 (более точно, к ее количественной версии, утверждению \ref{0-ra2m}).
Аналогично теорема~\ref{0-ne4} в размерности 4 сводится (в \S\ref{0ram4}) к теореме~\ref{rlt-cgslin} в размерности 3.

Результаты о неотъемлемых пересечениях дают естественное обобщение непланарности графов:
примеры \emph{двумерных} аналогов графов, не реализуемых ни в трехмерном, ни в четырехмерном пространствах.
Это разъяснено в замечании~\ref{r:hyper}.

Мы приводим упрощенное изложение, доступное неспециалистам в этой области.
Результаты формулируются в терминах систем точек.
Поэтому мы не используем понятия реализуемости гиперграфа (но упоминаем это понятие, поскольку оно важно для \emph{мотивировки}).
Для понимания большей части статьи достаточно знать основы стереометрии, а также знать или быть готовым изучить простейшие четырехмерные обобщения.
Мы уверены, что элементарное описание простых приложений топологических методов делает эти методы более доступными.
О сравнении с другими доказательствами см. замечание~\ref{s:compar}.

Описываемая связь между <<неотъемлемыми пересечениями>> и <<неотъемлемой зацепленностью>> не только дает простые доказательства классических результатов.
Она также подводит читателя к переднему краю науки, в частности, к доказательству обобщенной гипотезы Менгера (объясненному в \S\ref{s:prod}).
Изложенные результаты о неотъемлемой зацепленности являются отправной точкой \emph{(рамсеевской) теории неотъемлемых зацеплений}.
См. обзоры \cite{RA05, PS05, FMM+, Na20} и ссылки в них; о многомерных аналогах см. \cite{SS92, Sk03, KS20}.
Изложенные результаты о неотъемлемых пересечениях обобщаются до \emph{нереализуемости гиперграфов}.
См. обзоры \cite[\S4, \S5]{Sk06}, \cite[\S1]{MTW}, \cite[\S3.2]{Sk18}, 
и ссылки в них, а также недавние работы \cite{Pa15, Sk18o, AKM, PS20, Pa21, Me22, Sk23}.
Об аналогичных проблемах, относящихся к вложениям динамических систем см. \cite{LT14} и ссылки в этой работе.

История вопроса изложена в замечании~\ref{r:hist}.

\smallskip
\textbf{О независимости.}
Формально, замечания в дальнейшем  не используются и могут быть пропущены.
То же верно для \S\ref{s:stat4} и \S\ref{s:imrem}.
Параграфы \ref{0pro}, \ref{s:imrem} и \ref{s:prod} независимы друг от друга, так что их можно читать в любом порядке.
Ссылки на последующий текст, как и ссылки на другие работы, могут быть проигнорированы при первом чтении.

\subsection{Обобщения на высшие размерности и высшую кратность}\label{s:stat4}

Подмножество в $\R^d$ называется \emph{выпуклым}, если для любых двух точек этого подмножества отрезок, их соединяющий, также лежит в этом подмножестве.
\emph{Выпуклой оболочкой} множества $X \subset \R^d$ называется минимальное выпуклое множество, содержащее $X$.

\begin{theorem}\label{stat-il} Даны $d+3$ точки в $\R^d$.

Если $d=2k$, то среди них найдутся два непересекающиеся $(k+1)$-элементных подмножества, выпуклые оболочки которых пересекаются.

Если $d=2k-1$ и никакие $2k$ из данных точек не лежат в одной $(2k-2)$-мерной гиперплоскости, то найдется зацепленная пара $(k+1)$-элементных подмножеств с вершинами в данных точках (зацепленность означает, что выпуклая оболочка первого подмножества пересекает границу выпуклой оболочки второго подмножества в нечетном числе точек).
\end{theorem}

Эта теорема обобщает утверждение~\ref{0-ra2} (возьмите $d=2$), теорему~\ref{rlt-cgslin} (возьмите $d=3$) и теорему~\ref{0-ne4} (возьмите $d=4$).
Она принадлежит ван Кампену--Флоресу для четных $d>2$ \cite{vK32, Fl34} и Сегалу--Спешу--Ловасу--Шриверу--Танияме для нечетных $d>3$ \cite[Corollary 1.1]{LS98}, \cite{Ta00} (рассуждение про индексы в \cite[\S1]{SS92} допускает простое обобщение до теоремы~\ref{stat-il}; поэтому теорема~\ref{stat-il} для нечетных $d>3$ неявно доказана в \cite{SS92}).

Теорема~\ref{stat-il} доказывается индукцией по $d$.
База $d=1$ тривиальна.
Шаг индукции доказан в \S\ref{0pro} для $d=2,3,4$; доказательство для общего случая аналогично.

Аналог теоремы~\ref{stat-il} для $d+2$ точек

$\bullet$ не имеет смысла для $d=2k-1$, поскольку среди $2k+1$ точек в $\R^{2k-1}$ невозможно  взять два зацепленных подмножества;

$\bullet$ неверен для четных $d$ аналогично контрпримеру к соответствующему аналогу теоремы~\ref{0-ne4}.

Для нечетных $d$ имеется утверждение~\ref{0-ne3}.b о неотъемлемых пересечениях и его многомерный аналог.
Они слабее соответствующих теорем~\ref{rlt-cgslin} и~\ref{stat-il} о неотъемлемой зацепленности.
Другие результаты приведены в \cite[\S4]{Sk16}.


Сформулируем маломерный аналог вышеприведенных результатов для трехкратных неотъемлемых пересечений.

\begin{theorem}[\cite{Sa91g}]\label{t:vkfr} Для любых $11$ точек в пространстве найдутся три треугольника с вершинами в этих точках, имеющие общую точку, никакие два из которых не имеют общих вершин.
\end{theorem}


Удивительно, что известное доказательство такого элементарного результата использует  алгебраическую топологию.
Интересно найти элементарное доказательство.

Аналог теоремы~\ref{t:vkfr} для $10$ точек неверен. 
Действительно, в пространстве возьмем вершины тетраэдра и его центр, см. рис. \ref{f-full5}.
Либо возьмем каждую из этих $5$ точек с кратностью два, либо для каждой из них возьмем близкую к ней.

О многомерном $r$-кратном аналоге теорем \ref{stat-il} о пересечении и~\ref{t:vkfr} см.
\cite[Theorem 1.6]{Sk16}, \cite[гипотеза 3.1.4 и текст после нее]{Sk18}.
Об аналоге для $r$-кратной зацепленности см. \cite[\S5]{AKS}.

\subsection{Некоторые важные замечания}\label{s:imrem}

\begin{remark}[связь с гиперграфами]\label{r:hyper}
(a) Двумерными аналогами (конечных) графов являются \emph{$3$-однородные}, или \emph{двумерные} \emph{гиперграфы}, определяемые как системы трехэлементных подмножеств конечного множества.\footnote{В топологии такие объекты называются \emph{чистыми, или размерностно однородными}, двумерными \emph{симплициальными комплексами}.
Термин <<гиперграф>> знаком большинству математиков  и специалистам по компьютерным наукам.}
Для краткости, мы опускаем <<$3$-однородный>> и <<двумерный>>.
Например, \emph{полным гиперграфом} на $k$ вершинах называется семейство всех трехэлементных подмножеств $k$-элементного множества.
\emph{Реализуемость} (также называемая вложимостью) гиперграфа в $d$-мерном пространстве определяется аналогично реализуемости графа на плоскости: надо <<рисовать>> треугольники для каждого трехэлементного подмножества.
См. рис. \ref{f-full5}, \ref{f-mnre} и~\ref{f:k5i}; на последних рисунках подразбиение четырехугольников, аналогичное рис.~\ref{f-mnre} слева, не показано.
См. строгие определения, например, в \cite[\S3.2]{Sk18}, \cite[\S3]{Sk23'}.

Гиперграфы (и симплициальные комплексы) играют важную роль в математике.
Без них невозможно представить топологию и комбинаторику.
Кроме того, они применяются в компьютерных науках и биоинформатике, см. например \cite{PS11}.

Соображения <<малого шевеления>> (или <<общего положения>>) показывают, что любой граф реализуем в $\R^3$.
Несложное обобщение показывает, что любой гиперграф реализуем в $\R^5$.

Полный гиперграф на $6$ вершинах содержит <<конус над $K_5$>>, а потому не реализуем в $\R^3$ (утверждение~\ref{0-ne3}.a).
Уже в ранней истории топологии (1920-е гг.) математики пытались построить гиперграфы, не реализуемые в $\R^4$.
Эгберт ван Кампен и Алехандро Флорес в 1932--34 гг. доказали, что полный гиперграф на $7$ вершинах не реализуем в $\R^4$ (ср. теорему~\ref{0-ne4}).
Этот результат является и ранним приложением \emph{комбинаторной топологии} (сегодня она называется алгебраической топологией), и одним из первых результатов \emph{топологической комбинаторики} (также являющейся областью современных активных исследований).

(b) Реализациями (=вложениями) называются отображения без самопересечений.
Для топологической комбинаторики и дискретной геометрии интересно изучать отображения, самопересечения которых непусты (в отличие от вложений), но <<не слишком замысловаты>>.
Важным частным случаем является изучение отображений \emph{без тройных пересечений} и вообще отображений \emph{без $r$-кратных пересечений}, см. 
\S\ref{s:stat4} и обзоры \cite{Sk16}, \cite[\S3.3]{Sk18}, \cite[\S3]{Sk23'}.

(c) Мы представляем \emph{линейные} версии результатов.
О \emph{кусочно-линейной (PL)} и \emph{топологической} реализуемости гиперграфов см., например, \cite[\S3.2]{Sk18}, \cite[\S3]{Sk23'}, \cite[\S6 <<Реализуемость гиперграфов и комплексов>>]{Sk}.
Приводимые доказательства интересны тем, что они легко обобщаются до PL версий \cite{Sk03, Z13}, в отличие от доказательств из \cite{BM15, So12}.

Из PL версий <<количественных>> результатов (см. \S\ref{s:quant}) вытекают PL версии для \emph{почти-вложений} (см. PL случай теорем 1.4.1 и 3.1.6 из \cite{Sk18}, \cite{Sk23'}).
Последние влекут \emph{топологические} версии (см. объяснения в \cite[абзац после теоремы 1.4.1]{Sk18}).
\end{remark}

\begin{remark}[сравнение с другими доказательствами]\label{s:compar}
Теорема~\ref{stat-il} для четного $d$ (а потому и ее частные случаи --- утверждение \ref{0-ra2} и теорема \ref{0-ne4}) имеет другое простое доказательство, использующее \emph{число ван Кампена}, см. например \cite[\S1.4]{Sk18}, \cite[\S1.4]{Sk23'}, \cite[\S1.4, \S6 <<Реализуемость гиперграфов и комплексов>>]{Sk}.
Эти два доказательства предположительно являются простейшими из известных (<<доказательства из Книги>>); они
проходят для количественной, PL и топологической версий, ср. замечание \ref{r:hyper}.c.

Обычно теорема \ref{stat-il} для четного $d$ (точнее, топологическая версия ее и теоремы \ref{0-ne4j}) доказывается при помощи теоремы Борсука-Улама \cite[\S8]{Sk20}, \cite[\S5]{Ma03}.
В отличие от доказательства в данной статье (и от альтернативного  доказательства, использующего число ван Кампена), это требует некоторого знания алгебраической топологии.
Это знание не упрощает жизнь: известные доказательства теоремы Борсука-Улама (см. \cite{Ma03} и ссылки там) не проще, чем обсужденные выше прямые доказательства теоремы \ref{stat-il} для четного $d$.
(Теорема Борсука-Улама доказывается с использованием понятия \emph{степени} аналогично прямому доказательству теоремы~\ref{stat-il} для четного $d$, использующему \emph{число ван Кампена}.)

В \cite{BM15} представлено короткое алгебраическое доказательство теоремы \ref{stat-il} (и тем самым ее частных случаев).
Это доказательство в духе алгебраического доказательства \emph{теоремы Радона}\footnote{См., формулировку теоремы Радона, например, в \cite[\S3]{Sk23'}, \cite[\S1]{Sk16}.
О связях между теоремой Радона и теоремами \ref{rlt-cgslin}, \ref{0-ne4}, \ref{stat-il} см. \cite[\S3]{Sk23'}, \cite[\S4]{Sk16}.
Доказательство в \cite{BM15}, по-видимому, является прямой (т.~е., не использующей преобразования Гейла) версией доказательства теоремы~\ref{stat-il} для четного $d$ (это \cite[Theorem 5]{So12} для $k=1$).}.
\end{remark}

\begin{remark}[историческое]\label{r:hist}
Общие идеи <<понижения размерности>> или <<линка вершины>> просты и хорошо известны (см. замечание~\ref{s:lowering}).
О доказательствах теоремы Радона, основанных на этой идее, см. \cite{Pe72, Ko18, RRS}.
О приложении в компьютерной науке см. \cite[proof of 2.3.i]{DE94}.
Также известна связь между зацепленностью и пересечениями.\footnote{Например, коэффициент зацепления двух непересекающихся замкнутых ломаных в трехмерной сфере $\partial D^4$ равен алгебраическому числу пересечений двух двумерных дисков \emph{общего положения} в четырехмерном шаре $D^4$, натянутых на эти ломаные.
Индуктивное доказательство, включающее утверждение о зацепленности в нечетных размерностях и утверждение о пересечениях в четных размерностях, приведено в \cite[Whitney Lemma 5.12 and Theorem 5.16]{RS72}.}
Развитие этой идеи до связи между неотъемлемой зацепленностью и нереализуемостью нетривиально (ср. с разницей между утверждением \ref{0-ne3}.a и теоремой \ref{rlt-cgslin}).
Доказательства, открывшие и использовавшие эту связь, видимо, не были опубликованы

$\bullet$ до \cite{RST91, RST93, Sh03, RSSZ, Z13} --- для доказательства теоремы Конвея--Гордона--Закса \ref{rlt-cgslin} сведением неотъемлемой зацепленности к неотъемлемым пересечениям в низшей размерности;

$\bullet$ до \cite[Example 2, Lemmas 2 and 1']{Sk03}, \cite{RSSZ} --- для доказательств теоремы \ref{0-ne4} и гипотезы Менгера (см. \S\ref{s:menger}) сведением неотъемлемых пересечений к неотъемлемой зацепленности в низшей размерности.
\end{remark}

\section{Доказательства}\label{0pro}

Под \emph{$k$ точками в $\R^d$} (в этой статье в основном $d\le4$) подразумевается $k$ попарно различных точек.
Обозначим $[n]:=\{1,2,\ldots,n\}$.

\subsection{Плоскость: доказательство утверждения~\ref{0-ra2}}\label{0pla}


Утверждение~\ref{0-ra2} легко доказывается перебором количества вершин в выпуклой оболочке $5$ точек.
Чтобы продемонстрировать идею <<понижения размерности>> (см. замечание~\ref{s:lowering}) в наипростейшей ситуации, 
мы выведем утверждение~\ref{0-ra2} из утверждения~\ref{p:liline}. 

\begin{figure}[h]\centering
 \includegraphics{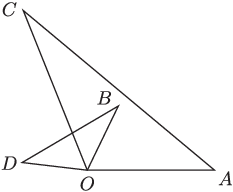} \qquad \qquad\includegraphics[scale=0.6]{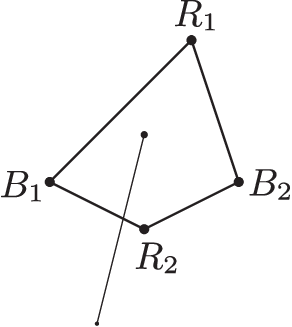}
\caption{К доказательству утверждений~\ref{0-ra2} (слева) 
и~\ref{0-k33}.b (справа).}
\label{f:plane}
\end{figure}

\begin{proof}[Доказательство утверждения \ref{0-ra2}]
Существует прямая, относительно которой ровно одна точка $O$ из данных $5$ лежит по одну сторону, а остальные точки $A,B,C,D$ "--- по другую.
Если для некоторых двух точек $X,Y\in\{A,B,C,D\}$ точка $X$ лежит на отрезке $OY$, то утверждение доказано.
В противном случае можно считать, что точки $A,B,C,D$ видны из точки $O$ именно в таком порядке, см. рис.~\ref{f:plane}.
Тогда получаем нужное по нижеприведенной лемме~\ref{0-eveng}.
\end{proof}

\begin{lemma}[о понижении размерности; см. рис.~\ref{f:plane} слева]\label{0-eveng}
Даны два треугольника на плоскости, имеющие единственную общую вершину.
Прямая $l$ отделяет ее от оснований треугольников.
Пересечения прямой $l$ с контурами треугольников чередуются вдоль прямой.
Тогда некоторые две стороны треугольников пересекаются, не имея общих вершин.
\end{lemma}

Эта лемма легко следует из утверждения~\ref{p:liline}.
Она сформулирована явно для удобства ее применения (здесь и в \S\ref{0ram4}) и для иллюстрации ее многомерного обобщения (лемма~\ref{0-even4}).

Следующие утверждения доказываются аналогично утверждению~\ref{0-ra2}.
Они используются для некоторых трехмерных результатов (утверждения~\ref{p:k4k4}.ab и теорема~\ref{rlt-sachs}).

\begin{proposition}\label{0-k33} (a) (рис.~\ref{k5} справа) Даны две тройки точек на плоскости.
Тогда найдутся два пересекающиеся отрезка, не имеющие общих вершин, каждый из которых соединяет точки из разных троек.

(b) (рис.~\ref{f:plane} справа) Даны четыре красных и две синих точки $B_1,B_2$ на плоскости, 
для которых любые два отрезка с разноцветными концами либо не пересекаются, либо пересекаются по общему концу.
Тогда найдутся две красные точки $R_1,R_2$ такие, что четырехугольник $R_1B_1R_2B_2$ несамопересекающийся, 
а две оставшиеся  красные точки лежат по разные стороны от этого четырехугольника.
(<<По разные стороны>> означает, что ломаная \emph{общего положения}, соединяющая две оставшиеся красные точки, 
пересекает контур четырехугольника в нечетном числе точек.)
\end{proposition}

\subsection{Ослабленные версии теоремы \ref{rlt-cgslin}}\label{s:weaker}

Сначала мы проиллюстрируем идею <<понижения размерности>> (см. замечание~\ref{s:lowering}) доказательством следующих ослабленных версий теоремы~\ref{rlt-cgslin}.

\begin{proposition}\label{0-ne3} Среди любых $6$ точек в пространстве найдутся

(a) $5$ точек $O,A,B,C,D$, для которых треугольники $OAB$ и $OCD$ имеют общую точку, отличную от $O$;

(b) непересекающиеся пара и тройка, для которых отрезок, соединяющий вершины пары, пересекает треугольник, образованный тройкой.
\end{proposition}

\begin{figure}[h]\centering
\begin{tikzpicture}[line cap=round,line join=round,>=triangle 45,x=0.32cm,y=0.32cm]
\clip(-4.93,-7.56) rectangle (15.43,4.12);
\fill[fill=black,fill opacity=0.1] (-3,-5) -- (0,1) -- (13,1) -- (10,-5) -- cycle;
\draw (4,3)-- (1,0);
\draw (4,3)-- (-1,-4);
\draw (4,3)-- (3,-2);
\draw (4,3)-- (7,-2);
\draw (4,3)-- (11,0);
\draw (2.4,-5)-- (2,-7);
\draw [dash pattern=on 2pt off 2pt] (3,-2)-- (2.4,-5);
\draw [dash pattern=on 2pt off 2pt] (7,-2)-- (8.76,-5);
\draw (8.76,-5)-- (9.85,-6.86);
\draw [dash pattern=on 2pt off 2pt] (1,0)-- (-2,-3);
\draw (-2,-3)-- (-3,-4);
\draw (-1.71,-5)-- (-2.4,-5.96);
\draw (12.24,-0.53)-- (13.38,-1.02);
\draw [dash pattern=on 2pt off 2pt] (-1.71,-5)-- (-1,-4);
\draw [dash pattern=on 2pt off 2pt] (11,0)-- (12.24,-0.53);
\draw (-3,-5)-- (0,1);
\draw (13,1)-- (10,-5);
\draw (10,-5)-- (-3,-5);
\draw (0,1)-- (1.66,1);
\draw (3,1)-- (3.34,1);
\draw (4,1)-- (5,1);
\draw (6,1)-- (8,1);
\draw (9.32,1)-- (13,1);
\fill [color=black] (1,0) circle (1.5pt);
\fill [color=black] (-1,-4) circle (1.5pt);
\fill [color=black] (7,-2) circle (1.5pt);
\fill [color=black] (11,0) circle (1.5pt);
\fill [color=black] (3,-2) circle (1.5pt);
\fill [color=black] (4,3) circle (1.5pt);
\draw[color=black] (2.83,3.24) node {$O$};
\fill [color=black] (2,-7) circle (1.5pt);
\draw[color=black] (3.27,-6.36) node {$A_5$};
\fill [color=black] (9.85,-6.86) circle (1.5pt);
\draw[color=black] (8.69,-6.32) node {$A_4$};
\fill [color=black] (-3,-4) circle (1.5pt);
\draw[color=black] (-3.37,-2.64) node {$A_2$};
\fill [color=black] (-2.4,-5.96) circle (1.5pt);
\draw[color=black] (-1.09,-5.8) node {$A_3$};
\fill [color=black] (13.38,-1.02) circle (1.5pt);
\draw[color=black] (13.05,-2.14) node {$A_1$};
\end{tikzpicture} 
\caption{К доказательствам утверждения \ref{0-ne3}.a и теоремы~\ref{rlt-cgslinm}.
Плоскость в $\R^3$ пересекает отрезки $OA_1,\ldots,OA_5$ по точкам $A_1',\ldots,A_5'$.}
\label{f:rlt-cgslin}
\end{figure}

\begin{proof}[Доказательство части (a)]
Cуществует плоскость, по одну сторону от которой лежит ровно одна точка $O$ из данных, а остальные 5 точек "--- по другую.
Рассмотрим пересечение этой плоскости с объединением треугольников $OAB$ для всех пар $A,B$ данных точек.
Теперь часть (a) следует из утверждения~\ref{0-ra2}.
\end{proof}

Часть (b) вытекает из (a) (и наоборот).


Рисунок~\ref{f-full5} показывает, что аналог каждой из частей (a,b) для $5$ точек неверен.

\subsection{<<Количественные>> версии}\label{s:quant}

Мы доказываем следующие усиленные <<количественные>> (алгебраические) версии результатов из \S\ref{0intr}.
(Вывод теоремы~\ref{rlt-cgslin} из утверждения~\ref{0-ra2}, а не из его <<количественной>> версии, содержит перебор, 
который трудно обобщить на высшие размерности.)

{\renewcommand{\thetheorem}{\ref{p:liline}$'$}\addtocounter{theorem}{-1}
\begin{proposition}[очевидное]\label{p:liline'}
Любые $4$ точки на прямой можно единственным (с точностью до порядка) способом разбить на две незацепленные пары.
\end{proposition}
}

{\renewcommand{\thetheorem}{\ref{0-ra2}$'$}\addtocounter{theorem}{-1}
\begin{proposition}\label{0-ra2m}
     Даны $5$ точек на плоскости, никакие $3$ из которых не лежат на одной прямой.
Рассмотрим все отрезки с вершинами в данных $5$ точках.
Тогда количество тех точек пересечения таких отрезков, которые отличны от вершин, нечетно.
\end{proposition}
}

Это легко доказывается перебором выпуклых оболочек данных $5$ точек, или вытекает из утверждения \ref{p:liline'} и следующей леммы.

{\renewcommand{\thetheorem}{\ref{0-eveng}$'$}\addtocounter{theorem}{-1}
\begin{lemma}\label{0-evengq}
(a) Даны два треугольника на плоскости, имеющие единственную общую вершину $O$.
    Никакие $3$ из их $5$ вершин не лежат на одной прямой.
    Прямая $l$ отделяет вершину $O$ от оснований треугольников.
    Тогда контуры треугольников пересекаются в четном числе точек тогда и только тогда, когда пересечения прямой $l$ с контурами треугольников чередуются вдоль прямой (т.е. когда пересечение одного треугольника с контуром другого содержит ровно один отрезок с концом в $O$).

(b) Даны $5$ точек $O,A_1,A_2,A_3,A_4$ на плоскости, никакие $3$ из которых не лежат на одной прямой.
Тогда количество из утверждения \ref{0-ra2m} равно сумме количеств точек пересечения внутренностей сторон  треугольников $OPQ$ и $ORS$, по всем неупорядоченным разбиениям точек $A_1,A_2,A_3,A_4$ на две неупорядоченные пары $P,Q$ и $R,S$.

(c) \cite[Proposition 7.5.a]{DGN+} Обозначим через $X={{[5]\choose2}\choose2}$ множество всех неупорядоченных пар двухэлементных подмножеств множества $[5]$.
Для каждого (из всех трех) неупорядоченного разбиения $\sigma\sqcup\tau=[4]$ на двухэлементные множества обозначим
    \[
        T_{\{\sigma,\tau\}}:=\bigl\{\{\alpha,\beta\}\in X\ :\  \alpha\subset\sigma\sqcup\{5\},\ \beta\subset\tau\sqcup\{5\}\bigr\}.
    \]
При этих обозначениях пара $\{\alpha,\beta\}\in X$ лежит в нечетном числе множеств $T_{\{\sigma,\tau\}}$ (т.е. в одном или в трех из этих трех множеств) тогда и только тогда, когда $\alpha\cap\beta=\emptyset$.\footnote{Иными словами, сумма по модулю $2$ (симметрическая разность) множеств $T_{\{\sigma,\tau\}}$ по всем таким разбиениям $\sigma\sqcup\tau=[4]$ равна $\bigl\{\{\alpha,\beta\}\in X\  :\ \alpha\cap\beta=\emptyset\bigr\}$.}
\end{lemma}
}

Часть (a) аналогична лемме \ref{0-eveng} (и тривиальна).
Часть (b) следует из (c).
Несложное доказательство части (c) оставляем читателю.

\begin{remark}\label{r:stronger}
Утверждение~\ref{0-ra2m} действительно сильнее утверждения~\ref{0-ra2}, поскольку достаточно доказать утверждение~\ref{0-ra2} при дополнительном условии, что никакие $3$ из $5$ точек не лежат на одной прямой,

(a)  во-первых, так как иначе утверждение~\ref{0-ra2} очевидно:
если точки $A,B,C$ среди данных пяти лежат на одной прямой, $B$ лежит между $A$ и $C$, а $D$ --- любая другая оставшаяся точка, то отрезки $AC$ и $BD$ пересекаются,

(b) во-вторых, так как можно немного пошевелить данные точки, чтобы никакие $3$ из $5$ пошевеленных точек не лежали на одной прямой и чтобы не добавилось точек пересечения отрезков с различными концами.

Аналогично рассуждению (b) получам, что нижеприведенные теоремы \ref{rlt-cgslinm} и \ref{0-ne4m} сильнее теорем~\ref{rlt-cgslin} и~\ref{0-ne4}.
\end{remark}

Треугольники, отличающиеся перестановкой вершин, считаются одинаковыми.

{\renewcommand{\thetheorem}{\ref{rlt-cgslin}$'$}\addtocounter{theorem}{-1}
\begin{theorem}[\cite{Sa81, CG83}]\label{rlt-cgslinm}
Даны $6$ точек в пространстве, никакие $4$ из которых не лежат в одной плоскости.
Тогда количество неупорядоченных пар зацепленных треугольников с вершинами в данных $6$ точках нечетно.
\end{theorem}
}

{\renewcommand{\thetheorem}{\ref{0-ne4}$'$}\addtocounter{theorem}{-1}
\begin{theorem}[\cite{vK32, Fl34}]\label{0-ne4m}
Даны $7$ точек в $\R^4$, никакие $5$ из которых не лежат в одной трехмерной гиперплоскости. 
Рассмотрим все треугольники с вершинами в данных $7$ точках. 
Тогда количество тех точек пересечения таких треугольников, которые не лежат на объединении ребер, нечетно.
\end{theorem}
}

Теорема \ref{stat-il} имеет аналогичную количественную версию.

Другие количественные версии приведены в п. \ref{s:unl}.
О контрпримерах к <<целочисленным версиям>> см.~\cite[Proposition~1.2 and Theorem~1.4]{KS20}.

\subsection{Трехмерное пространство: доказательство теоремы \ref{rlt-cgslinm}}\label{0ramcgs}


\begin{proof}[Доказательство теоремы \ref{rlt-cgslinm}]
Существует плоскость $\alpha$, по одну сторону от которой лежит ровно одна точка $O$ из данных, а остальные точки $A_1,\ldots,A_5$ "--- по другую (рис. \ref{f:rlt-cgslin}).
Возьмем 5 точек пересечения плоскости $\alpha$ с отрезками $OA_1,\ldots,OA_5$.
Так как никакие $4$ из данных $6$ точек не лежат в одной плоскости, то никакие $3$ из взятых $5$ точек не лежат на одной прямой.
При взгляде из точки $O$ на плоскость $\alpha$ получаем картину, аналогичную рис.~\ref{genpos}, посередине.

В пространстве \emph{отрезок $MN$ ниже отрезка $q$} (при взгляде из точки $O$), если отрезок $q$ пересекает внутренность треугольника $OMN$.   
Треугольники $OA_1A_2$ и $A_3A_4A_5$ зацеплены тогда и только тогда, когда отрезок $A_1A_2$ ниже \emph{ровно одной стороны} треугольника $A_3A_4A_5$.
Отрезок не может пересекать треугольник более, чем по двум точкам. 
Поэтому  <<ровно одной стороны>> можно заменить на <<нечетного числа сторон>>.

Тогда следующие числа имеют одинаковую четность:

$\bullet$ количество зацепленных неупорядоченных пар треугольников, образованных данными $6$ точками;

$\bullet$ количество отрезков $A_iA_j$, лежащих ниже нечетного количества сторон их <<дополнительных>> треугольников $A_kA_lA_m$, $\{i,j,k,l,m\}=\{1,2,3,4,5\}$;

$\bullet$ количество <<проходов>>, т.е. таких упорядоченных пар $(A_iA_j,A_kA_l)$ отрезков, у которых первый отрезок ниже второго;

$\bullet$ количество точек пересечения внутренностей отрезков с концами в 5 взятых точках в $\alpha$.
 
По утверждению~\ref{0-ra2m} последнее число нечетно.
Следовательно, первое число также нечетно.
\end{proof}


Следующая версия теоремы \ref{rlt-cgslin} аналогично сводится к утверждению \ref{0-k33}.b \cite{Z13}.
Эта версия используется для некоторого четырехмерного результата (теорема \ref{0-ne4pr}) в \S\ref{s:mnre4}.

\begin{figure}[h]\centering
\includegraphics[scale=0.6]{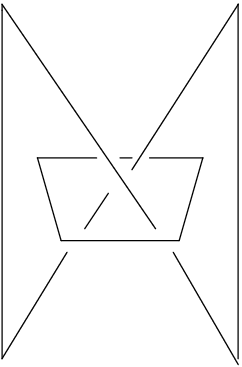}
\caption{Зацепление Уайтхеда, образованное пространственными четырехугольниками}
\label{f:whi}
\end{figure}

Возьмем два пространственных четырехугольника (т.~е. две замкнутые четырехзвенные ломаные) $ABCD$ и $A'B'C'D'$, никакие $4$ из $8$ вершин которых не лежат в одной плоскости.
Четырехугольники называются \emph{зацепленными по модулю $2$}, если количество точек пересечения четырехугольника $ABCD$ с объединением треугольников $A'B'C'$ и $A'D'C'$ нечетно.
(В отличие от случая треугольников, в пространстве существуют \emph{зацепленные} четырехугольники, не являющиеся при этом зацепленными по модулю $2$, см. рис.~\ref{f:whi}.
Ср. с определением в начале \S\ref{s:mnre4}.)

\begin{theorem}[\cite{Sa81}]\label{rlt-sachs}
Даны $4$ красные и $4$ синие точки в пространстве.
Никакие $4$ из этих $8$ точек не лежат в одной плоскости.
Тогда найдутся два зацепленных по модулю $2$ пространственных четырехугольника, состоящиx из отрезков с разноцветными концами.
\end{theorem}

\subsection{Четырехмерное пространство: доказательства теорем \ref{0-ne4} и~\ref{0-ne4m}}\label{0ram4}


\begin{remark}[наглядные замечания о четырехмерном пространстве]\label{s:intu4}
    (a) <<Типичное>> пересечение

    $\bullet$ двух отрезков на плоскости "--- либо пустое множество, либо точка (здесь под <<типичностью>> подразумевается, что никакие три из концов отрезков не лежат на одной прямой);

    $\bullet$ отрезка и треугольника в пространстве "--- либо пустое множество, либо точка;

    $\bullet$ двух треугольников в четырехмерном пространстве "--- либо пустое множество, либо точка.

    (b) Для любых двух точек

    $\bullet$ плоскости, отличных от некоторой точки $A$, найдется ломаная, соединяющая эти две точки и не проходящая через $A$;

    $\bullet$ пространства, не лежащих на некоторой прямой $l$, найдется ломаная, соединяющая эти две точки и не пересекающаяся с $l$;

    $\bullet$ четырехмерного пространства, не лежащих в некоторой двумерной плоскости $\alpha$, найдется ломаная, соединяющая эти две точки и не пересекающаяся с $\alpha$.
\end{remark}

Более продвинутая интуиция о четырехмерном пространстве не является здесь необходимой. 
Однако она может быть развита изучением, например, \cite[\S5.1 <<Как работать с четырехмерным пространством?>>]{Sk}, \cite[\S2.1]{Sk24}.

\begin{proof}[Доказательство теоремы \ref{0-ne4}]
Существует трехмерная гиперплоскость $\alpha$, по одну сторону от которой лежит ровно одна точка $O$ из данных, а остальные 6 точек "--- по другую (рис.  \ref{f:4space}).
Возьмем $6$ точек пересечения плоскости $\alpha$ с отрезками, соединяющими $O$ с остальными точками.
Можно считать, что никакие $5$ из данных $7$ точек не лежат в одной трехмерной гиперплоскости (аналогично замечанию \ref{r:stronger}.b).
Значит, никакие $4$ из взятых $6$ точек в $\alpha$ не лежат в одной плоскости.
Тогда по теореме~\ref{rlt-cgslin} имеются два зацепленных треугольника с вершинами во  взятых точках.
Теперь получаем нужное по нижеприведенной лемме~\ref{0-even4}.
\end{proof}

\begin{figure}[h]\centering
\begin{tikzpicture}[line cap=round,line join=round,>=triangle 45,x=0.32cm,y=0.32cm]
\clip(-4.93,-7.56) rectangle (15.43,4.12);
\fill[fill=black,fill opacity=0.1] (-3,-5) -- (0,1) -- (13,1) -- (10,-5) -- cycle;
\draw (1,0)-- (-1,-4);
\draw (-1,-4)-- (7,-2);
\draw (9,-4)-- (11,0);
\draw (3,-2)-- (4.51,-2.5);
\draw (9,-4)-- (5.01,-2.67);
\draw (7,-2)-- (5.54,-1.51);
\draw (4,3)-- (1,0);
\draw (4,3)-- (-1,-4);
\draw (4,3)-- (3,-2);
\draw (4,3)-- (7,-2);
\draw (4,3)-- (9,-4);
\draw (4,3)-- (11,0);
\draw (2.4,-5)-- (2,-7);
\draw [dash pattern=on 2pt off 2pt] (3,-2)-- (2.4,-5);
\draw [dash pattern=on 2pt off 2pt] (7,-2)-- (8.76,-5);
\draw (8.76,-5)-- (9.85,-6.86);
\draw [dash pattern=on 2pt off 2pt] (1,0)-- (-2,-3);
\draw (-2,-3)-- (-3,-4);
\draw (-1.71,-5)-- (-2.4,-5.96);
\draw (9.71,-5)-- (10.74,-6.43);
\draw (12.24,-0.53)-- (13.38,-1.02);
\draw [dash pattern=on 2pt off 2pt] (-1.71,-5)-- (-1,-4);
\draw [dash pattern=on 2pt off 2pt] (9,-4)-- (9.71,-5);
\draw [dash pattern=on 2pt off 2pt] (11,0)-- (12.24,-0.53);
\draw (-3,-5)-- (0,1);
\draw (13,1)-- (10,-5);
\draw (10,-5)-- (-3,-5);
\draw (0,1)-- (1.66,1);
\draw (3,1)-- (3.34,1);
\draw (4,1)-- (5,1);
\draw (6,1)-- (8,1);
\draw (9.32,1)-- (13,1);
\draw (3,-2)-- (6.25,-1.19);
\draw (7.22,-0.94)-- (11,0);
\draw (1,0)-- (1.5,-0.17);
\draw (1.89,-0.3)-- (2.85,-0.62);
\draw (3.57,-0.86)-- (5.02,-1.34);
\fill [color=black] (1,0) circle (1.5pt);
\fill [color=black] (-1,-4) circle (1.5pt);
\fill [color=black] (7,-2) circle (1.5pt);
\fill [color=black] (9,-4) circle (1.5pt);
\fill [color=black] (11,0) circle (1.5pt);
\fill [color=black] (3,-2) circle (1.5pt);
\fill [color=black] (4,3) circle (1.5pt);
\draw[color=black] (2.83,3.24) node {$O$};
\fill [color=black] (2,-7) circle (1.5pt);
\draw[color=black] (3.27,-6.36) node {$A_5$};
\fill [color=black] (9.85,-6.86) circle (1.5pt);
\draw[color=black] (8.69,-6.32) node {$A_4$};
\fill [color=black] (-3,-4) circle (1.5pt);
\draw[color=black] (-3.37,-2.64) node {$A_2$};
\fill [color=black] (-2.4,-5.96) circle (1.5pt);
\draw[color=black] (-1.09,-5.8) node {$A_3$};
\fill [color=black] (10.74,-6.43) circle (1.5pt);
\draw[color=black] (11.63,-5.48) node {$A_6$};
\fill [color=black] (13.38,-1.02) circle (1.5pt);
\draw[color=black] (13.05,-2.14) node {$A_1$};
\end{tikzpicture} 
\caption{К доказательству теоремы \ref{0-ne4}; гиперплоскость $\alpha$ в $\R^4$ изображена как плоскость в $\R^3$}
\label{f:4space}
\end{figure}

\begin{lemma}[о понижении размерности]\label{0-even4}
Даны два тетраэдра $\tau$ и $\tau'$ в $\R^4$, имеющие единственную общую вершину.
Трехмерная гиперплоскость $\alpha$ отделяет ее от оснований тетраэдров.
Контуры треугольников $\alpha\cap\tau$ и $\alpha\cap\tau'$ не пересекаются и зацеплены в гиперплоскости $\alpha$.
Тогда некоторые две грани тетраэдров пересекаются, не имея общих вершин.
\end{lemma}

Эта лемма не так очевидна, как ее маломерные аналоги (лемма~\ref{0-eveng} и аналогичный результат для треугольника и тетраэдра в пространстве), поскольку \emph{поверхность тетраэдра в $\R^4$ не разбивает пространство $\R^4$} (ср. с замечанием \ref{s:intu4}.b).

\begin{proof}[Доказательство леммы \ref{0-even4}]
Обозначим через $\gamma$ плоскость пересечения трехмерных гиперплоскостей, натянутых на тетраэдры.
Тогда $\alpha\cap\gamma$ "--- прямая пересечения плоскостей зацепленных треугольников $\alpha\cap\tau$ и $\alpha\cap\tau'$.
Значит, $\tau\cap\gamma$ и $\tau'\cap\gamma$ "--- треугольники с общей вершиной $O$, являющейся общей вершиной тетраэдров (рис. \ref{f:plane} слева).
Так как $\alpha$ отделяет $O$ от оснований тетраэдров, то $\alpha\cap\gamma$ отделяет $O$ от оснований этих треугольников.
Так как треугольники $\alpha\cap\tau$ и $\alpha\cap\tau'$ зацеплены, то точки пересечения прямой $\alpha\cap\gamma$ с контурами треугольников $\tau\cap\gamma$ и $\tau'\cap\gamma$ чередуются вдоль этой прямой \cite[задача 4.1.1]{Sk} \cite[\S1.1]{Sk24} (ср. рис. \ref{genpos} справа).
Следовательно, по лемме \ref{0-eveng} две стороны треугольников $\tau\cap\gamma$ и $\tau'\cap\gamma$ пересекаются, не имея общих вершин.
Не более одной из этих сторон содержит точку $O$.
Поэтому найденные стороны лежат в пересекающихся гранях тетраэдров, не имеющих общих вершин.
\end{proof}

Теорема~\ref{0-ne4m} вытекает из теоремы~\ref{rlt-cgslinm} и следующей леммы \ref{relation}.

{\renewcommand{\thetheorem}{\ref{0-even4}$'$}\addtocounter{theorem}{-1}
\begin{lemma}\label{relation}
(a) Даны два тетраэдра $\tau$ и $\tau'$ в $\R^4$, имеющие единственную общую вершину $O$.
Никакие $5$ из их $7$ вершин не лежат в одной трехмерной гиперплоскости.
Трехмерная гиперплоскость $\alpha$ отделяет $O$ от оснований тетраэдров.
Тогда поверхности тетраэдров $\tau$ и $\tau'$ пересекаются в четном числе точек тогда и только тогда, когда треугольники $\alpha\cap\tau$ и $\alpha\cap\tau'$ зацеплены в гиперплоскости $\alpha$ (т.~е. когда пересечение одного тетраэдра с поверхностью другого содержит ровно один отрезок с концом в $O$).

(b) Даны $7$ точек $O,A_1,\ldots,A_6\in\R^4$, никакие $5$ из которых не лежат в одной трехмерной гиперплоскости.
Тогда количество из теоремы~\ref{0-ne4m} равно сумме количеств точек пересечения внутренностей граней тетраэдров $O\Delta$ и $O\Delta'$ по всем неупорядоченным разбиениям точек $A_1,\ldots,A_6$ на две неупорядоченные тройки $\Delta$ и $\Delta'$.
(\emph{Внутренностью} треугольника называется дополнение до его контура.)

(c) Обозначим через $X={{[7]\choose3}\choose2}$ множество всех неупорядоченных пар трехэлементных подмножеств множества $[7]$.
Для каждого (из всех десяти) неупорядоченного разбиения $\sigma\sqcup\tau=[6]$ на трехэлементные множества обозначим
    \[
        T_{\{\sigma,\tau\}}:=\bigl\{\{\alpha,\beta\}\in X\ :\  \alpha\subset\sigma\sqcup\{7\},\ \beta\subset\tau\sqcup\{7\}\bigr\}.
    \]
При этих обозначениях пара $\{\alpha,\beta\}\in X$ лежит в нечетном числе множеств $T_{\{\sigma,\tau\}}$ тогда и только тогда, когда $\alpha\cap\beta=\emptyset$.\footnote{Иными словами, сумма по модулю $2$ множеств $T_{\{\sigma,\tau\}}$ по всем таким разбиениям равна
$\bigl\{\{\alpha,\beta\}\in X\  :\ \alpha\cap\beta=\emptyset\bigr\}$.}
\end{lemma}
}

Часть (a) аналогична лемме \ref{0-even4}.
Часть (b) следует из (c).
Несложное доказательство части (c) оставляем читателю.


Следующая четырехмерная версия утверждения \ref{0-k33}.a связана (аналогично теореме \ref{0-ne4}) с некоторым трехмерным результатом о неотъемлемой зацепленности \cite[Remark 2.5]{DS22}.

\begin{theorem}[\cite{vK32, Fl34}]\label{0-ne4j} Даны три тройки точек в $\R^4$.
Тогда найдутся два пересекающиеся треугольника без общих вершин,
вершины каждого из которых лежат в трех разных тройках.
\end{theorem}

\subsection{Свойства <<не зацепленности>>}\label{s:unl}


Здесь мы приводим количественные версии, утверждающие четность количества пересечений (или зацепленностей), ср. п. \ref{s:quant}.

\begin{proposition}\label{r:lint}
(2) Существуют 5 точек плоскости, никакие 3 из которых не лежат на одной прямой, и такие, что любой отрезок, соединяющий две из них, пересекает контур треугольника, образованного остальными тремя, в четном числе точек.
(Т.е. каждая пара точек <<не зацеплена>> с треугольником, образованным остальными тремя.)

(2') Даны 5 точек на плоскости, никакие 3 из которых не лежат на одной прямой.
Тогда количество отрезков, их соединяющих и пересекающих ровно в одной точке контур треугольника, образованного остальными тремя точками, четно.
\end{proposition}

В пространстве вместо свойств не зацепленности \ref{r:lint}.2,2' есть свойство зацепленности (теорема \ref{rlt-cgslin}) и следующие свойства не зацепленности.

\begin{proposition}\label{r:lint3}
(3) Существуют 6 точек в пространстве, никакие 4 из которых не лежат в одной плоскости, и такие, что любой отрезок, соединяющий две из них, пересекает поверхность тетраэдра, образованного остальными четырьмя, в четном числе точек.
(Т.е. каждая пара точек <<не зацеплена>> с тетраэдром, образованным остальными четырьмя.)

(3') Даны 6 точек в пространстве, никакие 4 из которых не лежат в одной плоскости.
Тогда количество отрезков, их соединяющих и пересекающих ровно в одной точке поверхность тетраэдра, образованного остальными четырьмя точками, четно.
\end{proposition}


\begin{proof}[Набросок доказательства]
(3) Возьмем точки, близкие к вершинам правильного октаэдра, или точки, близкие к вершинам треугольной призмы, или  точки на кривой моментов.

(3') Любые два треугольника, натянутые на две непересекающиеся тройки данных точек, либо не пересекаются, либо пересекаются по отрезку (не вырожденному в точку).
Концов таких отрезков четное число.
Концы таких отрезков --- в точности точки пересечения отрезков, соединяющих данные точки, с поверхностями <<дополнительных>> тетраэдров.

{\it Другое доказательство.}
Имеем
$$\sum_{\{A,B\}} |AB\cap\partial T_{AB}|_2 = \sum_{\{A,B\}} (|A\cap T_{AB}|_2+|B\cap T_{AB}|_2) =
\sum_{A}\sum_{B\ne A} |A\cap T_{AB}|_2=0.$$
Здесь

$\bullet$ $|S|_2$ --- четность количества элементов в конечном множестве $S$;

$\bullet$ $T_{AB}$ --- тетраэдр, образованный четырьмя данными точками, отличными от $A,B$;

$\bullet$ последнее равенство верно, поскольку для любой точки $A$ множество $\{T_{AB}\}_{B\ne A}$ является \emph{3-циклом}, ср.  \cite[замечание 1.3.6b]{Sk23}. 
\end{proof}

Утверждения \ref{0-ra2}, \ref{0-ra2m}, \ref{0-ne3}.b и \ref{r:lint3}.3' показывают, что при переходе от размерности 2 к размерности 3 свойство существования пересечения сохраняется, а четность числа пересечений меняется.
Трехмерные аналоги утверждений \ref{0-ra2}, \ref{0-ra2m} имеют более сильную форму: теоремы \ref{rlt-cgslin} и \ref{rlt-cgslinm}.

\begin{proposition}\label{r:lint4}
(4-3) Существуют 7 точек в $\R^4$, никакие 5 из которых не лежат в одной трехмерной гиперплоскости, и такие, что любой треугольник, образованный тремя из них, пересекает поверхность тетраэдра, образованного остальными четырьмя точками, в четном числе точек. \cite{St24}
(Т.е. никакой треугольник, образованный тремя из этих точек, <<не зацеплен>> с тетраэдром, образованным остальными четырьмя.)

(4'-3) Даны 7 точек в $\R^4$, никакие 5 из которых не лежат в одной трехмерной гиперплоскости.
Тогда количество треугольников, образованных тремя из них и пересекающих ровно в одной точке тетраэдр, образованный остальными четырьмя точками, четно.

(4-2) (гипотеза) Существуют 7 точек в $\R^4$, никакие 5 из которых не лежат в одной трехмерной гиперплоскости, и такие, что любой отрезок, соединяющий две из них, пересекает поверхность четырехмерного симплекса, образованного остальными пятью точками, в четном числе точек.
(Т.е. никакая пара точек <<не зацеплена>> с четырехмерным симплексом, образованным остальными пятью.)

(4'-2) Даны 7 точек в $\R^4$, никакие 5 из которых не лежат в одной трехмерной гиперплоскости.
Тогда количество точек пересечения отрезков, их соединяющих, с трехмерными поверхностями четырехмерных симплексов, образованных остальными пятью точками, четно.
\end{proposition}

\begin{proof}[Набросок доказательства]
(4-3) Возьмем точки на кривой моментов \cite{St24}.

(4'-3) Пометим все точки пересечения пар треугольников, натянутых на две непересекающиеся тройки данных точек (это точки из теоремы \ref{0-ne4m}).
Тогда искомое количество равно удвоенному числу помеченных точек.
Следовательно, это количество четно.

(4-2) Попробуйте взять точки на кривой моментов.

(4'-2) Аналогично другому доказательству утверждения \ref{r:lint3}.3'.
\end{proof}

\begin{conjecture}\label{r:lintd}
(d-k) Для любого $d\ne2k-1$ существуют $d+3$ точек в $\R^d$, никакие $d+1$ из которых не лежат в одной $(d-1)$-мерной гиперплоскости, и такие, что любой $k$-симплекс, натянутый на $k+1$ из них, пересекает поверхность $(d+1-k)$-симплекса, натянутого на остальные $d+2-k$ точек, в четном числе точек.

{\it Указание.} Попробуйте взять точки на кривой моментов.

(d'-k) Даны $d+3$ точек в $\R^d$, никакие $d+1$ из которых не лежат в одной $(d-1)$-мерной гиперплоскости.
Тогда при $d\ne2k-1$ количество $k$-симплексов, натянутых на $k+1$ из них и пересекающих ровно в одной точке поверхность $(d+1-k)$-симплекса, натянутого на остальные $d+2-k$ точек, четно.

{\it Указание.} Индукция по $k$; аналогично другому доказательству утверждения \ref{r:lint3}.3'.
\end{conjecture}

\section{Реализуемость произведений и гипотеза Менгера}\label{s:prod}

\subsection{Гипотеза Менгера}\label{s:menger}

\emph{(Декартовым) произведением} $F\times F'$ двух фигур $F,F'$ в $\R^3$ называется множество точек $(x,y,z,x',y',z')\in\R^6$ таких, что $(x,y,z)\in F$ и $(x',y',z')\in F'$.

\begin{figure}[h]
\centerline{\includegraphics[scale=1.05]{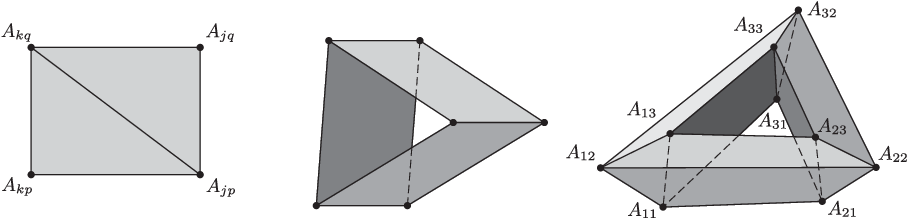}}
\caption{Реализации произведений: $K_2\times K_2$ (слева), $K_2\times K_3$ (в центре), $K_3\times K_3$ (справа)}
\label{f-mnre}
\end{figure}

\begin{figure}[h]
\includegraphics[width=8cm]{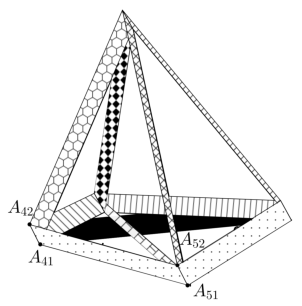}
\caption{Реализация произведения $K_5\times K_2$}
\label{f:k5i}
\end{figure}

Примеры реализации произведений приведены на рис.~\ref{f-mnre} и~\ref{f:k5i}.
Определение реализации см. например в \cite[\S3.2]{Sk18}, \cite[\S3]{Sk23'}, \cite[\S6 <<Реализуемость гиперграфов и комплексов>>]{Sk}.
Карл Менгер в 1929 году выдвинул гипотезу, что квадрат непланарного графа не реализуем в $\R^4$ \cite{Me29}
(ср. с теоремой \ref{0-ne4pr}).
Она была доказана только в 1978 году Брайаном Уммелем \cite{Um78}, использовавшим продвинутую алгебраическую топологию.
Простое же доказательство, использующее \emph{понижение размерности}, было получено в 2003 году Михаилом Скопенковым \cite{Sk03}, см. ниже.
Его рассуждение доказывает обобщенную гипотезу Менгера (<<$k$-ая степень непланарного графа не реализуема в $\R^{2k}$>>).
Оно даже дает простую формулу для минимального $d$ такого, что данное произведение нескольких графов реализуемо в $\R^d$ \cite{Sk03}.

Комбинаторной версией произведения является \emph{произведение двух графов} (не обязательно планарных).
Это произведение может быть рассмотрено (хотя и не канонически) как гиперграф.
См. замечание \ref{r:hyper}.a; ср. с \cite[\S3-\S6]{MNS}.

Из нижеприведенных утверждений \ref{p:k4k4}.ab вытекает, что ни $K_4\times K_4$, ни $K_5\times K_3$ не реализуемы (линейно) в $\R^3$.
Доказательство утверждения \ref{p:k4k4}.a показывает, что $K_{3,1}\times K_{3,1}$ не реализуемо в $\R^3$ (но это произведение реализуемо в $\R^4$).

\smallskip
\textbf{Замечание.}
Доказательство гипотезы Менгера, использующее число ван Кампена или теорему Борсука-Улама (см. замечание \ref{s:compar}), неизвестно.
Доказательство гипотезы Менгера в \cite{Um78} проходит для топологической версии, но сложно.
Более простое доказательство в \cite{Sk03} использует для топологической версии непростую  \emph{аппроксимационную теорему Брайанта}.
Доказательство топологической версии можно получить, придумав количественную PL версию гипотезы Менгера (т.~е., улучшив PL версию теоремы \ref{0-ne4pr} аналогично \S\ref{s:quant}, см. задачу \ref{c:menquan}).

\subsection{Реализуемость произведений}\label{0statpr}

Приведем формализацию идеи изображения в $\R^3$ или в $\R^4$ произведения $K_m\times K_n$.
Пусть имеется $mn$ точек $A_{jp}$, где $j\in[m]$ и $p\in[n]$,  в $\R^3$ или в $\R^4$.
Для чисел $j,k\in[m]$, $j<k$, и $p,q\in[n]$, $p<q$, обозначим
$$jk\times pq:=\{A_{jp}A_{kq}A_{jq},A_{jp}A_{kq}A_{kp}\}.$$
Это набор двух треугольников, имеющих общую сторону (рис.~\ref{f-mnre} слева).
Их объединение может быть, но не обязано быть, плоским четырехугольником.
Назовем \textbf{$(m,n)$-произведением} объединение наборов
\[
    jk\times pq,\quad\text{по всем}\quad 1\leqslant j<k\leqslant m,\quad 1\leqslant p<q\leqslant n.
\]
(Всего имеется $mn(m-1)(n-1)/2$ таких треугольников.)
\textbf{Телом} $(m,n)$-произведения называется объединение его треугольников.
Оно является многогранным, возможно, самопересекающимся

$\bullet$ квадратом, если $m=n=2$ (рис. \ref{f-mnre} слева);

$\bullet$ цилиндром (точнее, боковой поверхностью цилиндра), если $m=3$ и $n=2$ (рис. \ref{f-mnre} посередине);

$\bullet$ тором, если $m=n=3$ (рис. \ref{f-mnre} справа).

\begin{proposition}\label{p:k4k4} В пространстве в любом

(a) $(4,4)$-произведении; \quad (b) $(3,5)$-произведении

найдутся пересекающиеся треугольник и сторона треугольника, не имеющие общих вершин.
\end{proposition}

Эти фольклорные результаты сводятся в \S\ref{s:mnonre} к PL аналогам утверждений \ref{0-k33}.ab.

\begin{theorem}[о квадрате; \cite{Um78, Sk03}]\label{0-ne4pr}
В любом $(5,5)$-произведении в $\R^4$ найдутся два пересекающиеся треугольника, не имеющие общих вершин.
\end{theorem}

Теорема о квадрате \ref{0-ne4pr} сводится в \S\ref{s:mnre4} к версии теоремы \ref{rlt-sachs}.

\begin{example}\label{stat-epr} Аналоги утверждения~\ref{p:k4k4} и теоремы~\ref{0-ne4pr} неверны для

    (a) $(2,n)$-произведения в пространстве для каждого $n$ (для $n \leqslant 4$ это очевидно; для $n=5$ см. рис.~\ref{f:k5i}: вершины $(2,5)$-произведения являются вершинами параллелограммов; для $n \geqslant 6$ конструкция аналогична, см. \S\ref{s:mnre}; ср. \cite[Theorem 1.5]{RSS95'});

    (b) $(3,n)$-произведения в пространстве для каждого $n \leqslant 4$ (для $n \leqslant 3$ это очевидно, см. рис.~\ref{f-mnre} справа; для $n=4$ конструкция аналогична, см. \S\ref{s:mnre});

     (c) $(4,n)$-произведения в $\R^4$ для каждого $n$ (см. \S\ref{s:mnre}).
\end{example}

\begin{problem}[количественная проблема Менгера]\label{c:menquan}
Обозначим
\[
\t{K_5^2}:=\left\{\{(X,Y),(X',Y')\}\ :
    \ X,Y,X',Y'\in{[5]\choose2},\text{ либо } X\cap X'=\emptyset, \text{ либо }Y\cap Y'=\emptyset\right\}.
\]
Найдите подмножество $M\subset\t{K_5^2}$ такое, что
для любого \emph{кусочно-линейного отображения} $f\colon K_5\times K_5\to\R^4$ нечетно количество пар $\{(X,Y),(X',Y')\}\in M$, для которых количество точек в пересечении $f(X\times Y)\cap f(X'\times Y')$ нечетно, где под $X,Y,X',Y'$ понимаются ребра графа $K_5$.
\end{problem}

\subsection{Реализация произведений}\label{s:mnre}


\begin{proof}[Набросок доказательства примера \ref{stat-epr}.a]
Возьмем точки $A_{11},\ldots,A_{1n}\in\R^3$, никакие четыре из которых не лежат в одной плоскости.
Возьмем вектор $v$, не параллельный никакой плоскости, проходящей через три из этих точек.
Для каждого $p\in[n]$ возьмем точку $A_{2p}$, для которой $\vec{A_{1p}A_{2p}}=v$.
Если $v$ достаточно мал, то точки $A_{jp}$, $j\in\{1,2\}$, $p\in[n]$, "--- искомые, т.е. в образованном ими $(2,n)$-произведении не найдется пересекающихся треугольника и стороны треугольника, не имеющих общих вершин.

Действительно, $12\times pq$ "--- параллелограмм для любых $p\neq q$.
Поскольку никакие $4$ из точек $A_{11},\ldots,A_{1n}$ не лежат в одной плоскости, и по выбору вектора $v$, для любых попарно различных $p,q,r,s$ отрезки $A_{1p}A_{1q}$ и $A_{1r}A_{1s}$ не пересекаются.
Поскольку $v$ достаточно мал, предыдущее верно и с заменой $1$ на $2$.
Тогда любые два параллелограмма $12\times pq$ и $12\times rs$ (точнее, их выпуклые оболочки), не имеющие общей стороны, не пересекаются.
Теперь можно проверить, что точки $A_{jp}$ "--- искомые.
\end{proof}

\begin{proof}[Набросок доказательства примера \ref{stat-epr}.b]
    Обозначим
    $$A_{11}=(1,0,1),\quad A_{12}=(-1,0,1),\quad A_{13}=(0,0,2),\quad A_{14}=(0,0,3).$$
    Обозначим через $f\colon\R^3\to\R^3$ вращение на $\frac{2\pi}{3}$ вокруг оси $x$.
    Положим $A_{2p}=f(A_{1p})$ и $A_{3p}=f(f(A_{1p}))$ для всех $p\in[4]$ (рис.~\ref{f-mnre} справа).
    Тогда точки $A_{jp}$, $j\in[3]$, $p\in[4]$, "--- искомые.

Действительно, $jk\times pq$ "--- параллелограмм для любых $j\ne k, p\ne q$.
Поскольку любые два отрезка, соединяющие точки $A_{1p}$, либо не пересекаются, либо пересекаются по общей вершине, то любые два из этих параллелограммов, не имеющие общих сторон, не пересекаются.
Теперь можно проверить, что точки $A_{jp}$ "--- искомые.
\end{proof}

\begin{figure}[h]
\includegraphics[width=6cm]{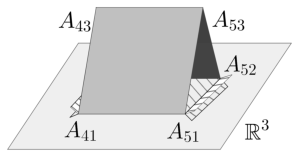}\qquad\includegraphics[width=6cm]{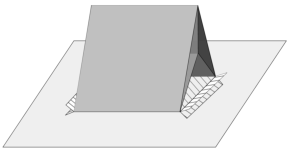}
\caption{Слева: к реализации произведения $K_3\times K_5$ в $\R^4$.
\newline
Справа: к реализации произведения $K_4\times K_5$ в $\R^4$.}
\label{f-k5xkn}
\end{figure}

\begin{proof}[Набросок доказательства ослабленной версии примера \ref{stat-epr}.c: $(3,5)$-произведение в $\R^4$]
    Возьмем трехмерную гиперплоскость в $\R^4$ (показанную на рис.~\ref{f-k5xkn} слева как плоскость в $\R^3$).
    В ней возьмем $10$ точек $A_{jp}$, где $j\in[5]$, $p\in\{1,2\}$, показанных на рис.~\ref{f:k5i}.
    Возьмем вектор $v$, не параллельный гиперплоскости.
    Положим $A_{j3}:=A_{j1}+v$.
    (На рис.~\ref{f-k5xkn} слева мы видим боковую поверхность призматоида $A_{41}A_{42}A_{43}A_{51}A_{52}A_{53}$.)
    Тогда точки $A_{jp}$, $j\in[5]$, $p\in[3]$, "--- искомые, т.е. в образованном ими $(3,5)$-произведении не найдется двух пересекающихся треугольников, не имеющих общих вершин.
\end{proof}

\begin{proof}[Набросок доказательства примера \ref{stat-epr}.c]
Возьмем точки $A_{jp}\in \R^3\subset \R^4$, $j\in\{1,2\}$, $p\in[n]$ из доказательства примера \ref{stat-epr}.a.
Тогда $\overrightarrow{A_{1p}A_{1q}}=\overrightarrow{A_{2p}A_{2q}}$ для всех $p\ne q$.
Возьмем неколлинеарные векторы $v_3,v_4\in\R^4$, не параллельные гиперплоскости $\R^3 \subset \R^4$.
Обозначим $A_{jp}:=A_{1p}+v_j$, $j\in\{3,4\}$.
Можно считать, что $v_3,v_4$ взяты так, что $A_{41}$ лежит внутри треугольника $A_{11}A_{21}A_{31}$.
См. рис.~\ref{f-k5xkn} справа.
Тогда точки $A_{jp}$, $j\in[4]$, $p\in[n]$, "--- искомые.
\end{proof}

\subsection{Нереализуемость произведений в трехмерном пространстве}\label{s:mnonre}

\begin{figure}[h]\centering
\includegraphics{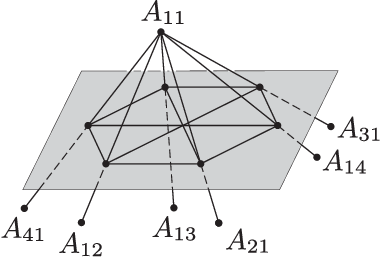}\qquad\begin{tikzpicture}[line cap=round,line join=round,>=triangle 45,x=0.32cm,y=0.32cm]
\clip(-4.91,-7.4) rectangle (15.45,4.56);
\fill[fill=black,fill opacity=0.1] (-3,-5) -- (0,1) -- (13,1) -- (10,-5) -- cycle;
\draw (-0.89,-3.76)-- (2.51,-3.71);
\draw (8.81,-3.6)-- (2.36,-3.71);
\draw (4,3)-- (0.89,-0.1);
\draw (4,3)-- (2.85,-1.32);
\draw (4,3)-- (-0.89,-3.76);
\draw (4,3)-- (10.95,0.08);
\draw (4,3)-- (8.81,-3.6);
\draw (4,3)-- (6.35,-2.52);
\draw (-1.79,-5)-- (-2.49,-5.96);
\draw [dash pattern=on 2pt off 2pt] (-0.89,-3.76)-- (-1.79,-5);
\draw [dash pattern=on 2pt off 2pt] (0.89,-0.1)-- (-2,-3);
\draw (-2,-3)-- (-3,-4);
\draw (1.87,-5)-- (1.39,-6.8);
\draw (9.83,-5)-- (10.48,-5.9);
\draw [dash pattern=on 2pt off 2pt] (1.87,-5)-- (2.85,-1.32);
\draw [dash pattern=on 2pt off 2pt] (8.81,-3.6)-- (9.83,-5);
\draw (-3,-5)-- (0,1);
\draw (13,1)-- (10,-5);
\draw (10,-5)-- (-3,-5);
\draw (0,1)-- (1.66,1);
\draw (3,1)-- (3.34,1);
\draw (4,1)-- (5,1);
\draw (6,1)-- (8,1);
\draw (9.32,1)-- (13,1);
\draw (2.85,-1.32)-- (-0.89,-3.76);
\draw (8.81,-3.6)-- (10.95,0.08);
\draw [dash pattern=on 2pt off 2pt] (6.35,-2.52)-- (7.4,-5);
\draw (7.4,-5)-- (8.15,-6.76);
\draw [dash pattern=on 2pt off 2pt] (10.95,0.08)-- (12.24,-0.51);
\draw (12.24,-0.51)-- (13.25,-0.96);
\draw (2.85,-1.32)-- (5.44,-0.88);
\draw (5.91,-0.8)-- (6.42,-0.71);
\draw (7.07,-0.6)-- (10.95,0.08);
\fill [color=black] (0.89,-0.1) circle (1.5pt);
\fill [color=black] (2.85,-1.32) circle (1.5pt);
\draw [color=black] (10.95,0.08) ++(-1.5pt,0 pt) -- ++(1.5pt,1.5pt)--++(1.5pt,-1.5pt)--++(-1.5pt,-1.5pt)--++(-1.5pt,1.5pt);
\fill [color=black] (8.81,-3.6) circle (1.5pt);
\fill [color=black] (6.35,-2.52) circle (1.5pt);
\draw [color=black] (-0.89,-3.76) ++(-1.5pt,0 pt) -- ++(1.5pt,1.5pt)--++(1.5pt,-1.5pt)--++(-1.5pt,-1.5pt)--++(-1.5pt,1.5pt);
\fill [color=black] (4,3) circle (1.5pt);
\draw[color=black] (3.07,3.76) node {$A_{11}$};
\fill [color=black] (-2.49,-5.96) circle (1.5pt);
\draw[color=black] (-1.13,-6.4) node {$A_{12}$};
\fill [color=black] (13.25,-0.96) circle (1.5pt);
\draw[color=black] (13.39,-1.94) node {$A_{13}$};
\fill [color=black] (-3,-4) circle (1.5pt);
\draw[color=black] (-3.21,-2.64) node {$A_{41}$};
\fill [color=black] (1.39,-6.8) circle (1.5pt);
\draw[color=black] (2.85,-6.64) node {$A_{21}$};
\fill [color=black] (10.48,-5.9) circle (1.5pt);
\draw[color=black] (11.53,-4.96) node {$A_{31}$};
\fill [color=black] (8.15,-6.76) circle (1.5pt);
\draw[color=black] (7.13,-6.56) node {$A_{51}$};
\end{tikzpicture}
\caption{К доказательствам утверждений \ref{p:k4k4}.a (слева) и \ref{p:k4k4}.b (справа) для $j=p=1$ [справа заменить $A_{ij}$ на $A_{ji}$ для всех $i\ne j$]}
\label{f:3space}
\end{figure}

\begin{proof}[Доказательство утверждения \ref{p:k4k4}.a]
(Аналогично утверждению~\ref{0-ne3}.)
Существует плоскость $\alpha$, по одну сторону от которой лежит ровно одна точка $A_{jp}$ из данных, а остальные 15 точек "--- по другую (рис.~\ref{f:3space} слева).
Точку пересечения плоскости $\alpha$

$\bullet$ с отрезком $A_{jp}A_{kp}$ покрасим в синий для каждого $k\in[4]-\{j\}$;

$\bullet$ с отрезком $A_{jp}A_{jq}$ покрасим в красный для каждого $q\in[4]-\{p\}$.

Пересечение произведения $jk\times pq$ с $\alpha$ назовем \emph{дугой}.
Тогда дуги имеют разноцветные концы.
(Пересечение плоскости $\alpha$ с телом данного $(4,4)$-произведения является PL изображением, возможно с  самопересечениями, графа $K_{2,4}$ в $\alpha$, т.е. образом \emph{PL отображения}  $K_{3,3}\to\alpha$.)

Значит, по PL аналогу\footnote{Этот аналог и доказывается аналогично, см. \cite[замечание 1.4.4b]{Sk23'}.
Можно взять $A_{11}$ вместо $A_{jp}$, маленькую сферу с центром в $A_{11}$ вместо $\alpha$, и сферический аналог утверждения~\ref{0-k33}.a вместо кусочно-линейного.
Ниже аналогичные замечания в доказательствах утверждения \ref{p:k4k4}.b и теоремы \ref{0-ne4pr} пропускаются.}
 утверждения~\ref{0-k33}.a
найдутся две пересекающиеся дуги, не имеющие общих концов.
Тогда найдутся $k,k'\in[4]-\{j\}$ и $q,q'\in[4]-\{p\}$ такие, что $k\ne k'$, $q\ne q'$, и некоторые два треугольника --- один из произведения $jk\times pq$, а другой из произведения $jk'\times pq'$ --- имеют общую точку, отличную от общей вершины $A_{jp}$ этих треугольников.
Значит, один из этих треугольников пересекает сторону другого треугольника, не проходящую через $A_{jp}$ и потому не имеющую с первым треугольником общих вершин.
\end{proof}

Для данных $9$ точек $A_{jp}$, $j\in\{u,v,w\}$, $p\in\{u',v',w'\}$, в $\R^3$ или в $\R^4$, обозначим через $uvw\times u'v'w'$ соответствующее $(3,3)$-произведение (рис.~\ref{f-mnre} справа; в отличие от рисунка, $(3,3)$-произведение может <<иметь самопересечения>>).

\begin{proof}[Доказательство утверждения \ref{p:k4k4}.b]
Можно считать, что никакие 4 из данных $15$ точек $A_{jp}$, $j\in[3]$, $p\in[5]$, не лежат в одной плоскости
(аналогично замечанию~\ref{r:stronger}.b).
Существует плоскость $\alpha$, по одну сторону от которой лежит ровно одна точка $A_{jp}$ из данных, а остальные 14 точек "--- по другую (рис.~\ref{f:3space} справа).
Точку пересечения плоскости $\alpha$

$\bullet$ с отрезком $A_{jp}A_{kp}$ покрасим в синий для каждого $k\in[3]-\{j\}$;

$\bullet$ с отрезком  $A_{jp}A_{jq}$ покрасим в красный для каждого $q\in[5]-\{p\}$.

Пересечение произведения $jk\times pq$ с $\alpha$ назовем \emph{дугой}.
Тогда дуги имеют разноцветные концы.
(Пересечение плоскости $\alpha$ с телом данного $(3,5)$-произведения является PL изображением, возможно с  самопересечениями, графа $K_{2,4}$ в $\alpha$.)

Аналогично последнему абзацу доказательства утверждения \ref{p:k4k4}.a либо

(1) в данном $(3,5)$-произведении найдутся пересекающиеся треугольник и сторона треугольника, не имеющие общих вершин, либо

(2) любые две дуги пересекаются только в их общей вершине (если она есть).

Во втором случае обозначим синие точки через $K$ и $L$.
По PL аналогу утверждения~\ref{0-k33}.b найдутся две красные точки $Q,R$ такие, что оставшиеся две красные точки $S,T$ лежат в $\alpha$ по разные стороны от замкнутой ломаной $\gamma$, являющейся объединением дуг $QK,KR,RL,LQ$.
Возьмем $k,l\in[3]-\{j\}$ и $q,r,s,t\in[5]-\{p\}$ такие, что точки $K,L$ и $Q,R,S,T$ принадлежат отрезкам, соединяющим $A_{jp}$ с $A_{kp},A_{lp}$ и $A_{jq},A_{jr},A_{js},A_{jt}$ соответственно.
Тогда плоскость $\alpha$ пересекает

$\bullet$ контур треугольника $j\times pst:=A_{jp}A_{js}A_{jt}$ по точкам $S$ и $T$ (заметим, что этот треугольник не лежит в данном $(3,5)$-произведении);

$\bullet$ тело $(3,3)$-подпроизведения $jkl\times pqr$ (содержащегося в данном $(3,5)$-произведении) по ломаной $\gamma$.

Если $\gamma$ имеет самопересечения, то получаем свойство (1).
Если не имеет, то так как точки $S,T$ лежат в $\alpha$ по разные стороны от $\gamma$,
по нижеприведенной лемме \ref{0-even3} получаем свойство (1).
\end{proof}

\begin{lemma}[о понижении размерности]\label{0-even3}
Пусть в пространстве контур $\partial\Delta$ треугольника $\Delta$ и $(3,3)$-произведение $\tau$ имеют единственную общую вершину $O$.
Никакие $4$ из вершин треугольника $\Delta$ и $(3,3)$-произведения $\tau$ не лежат в одной плоскости.
Плоскость $\alpha$ отделяет $O$ от основания треугольника $\Delta$, и от оставшихся $8$ точек $(3,3)$-произведения $\tau$.
Плоскость $\alpha$ пересекает

$\bullet$ тело $|\tau|$ по замкнутой ломаной, не имеющей самопересечений;

$\bullet$ $\partial\Delta$ по двум точкам $S,T$, лежащим по разные стороны от ломаной.

Тогда в $\tau$ и в $\partial\Delta$ найдутся пересекающиеся треугольник и сторона треугольника, не имеющие общих вершин.
\end{lemma}

\begin{proof} Обозначим через $\left<\Delta\right>$ плоскость треугольника $\Delta$.
Тогда $\alpha\cap\left<\Delta\right>$ "--- прямая.
Пересечение $\alpha\cap\Delta$ "--- отрезок $ST$.
Точки $S,T$ лежат в $\alpha$ по разные стороны от замкнутой ломаной $\alpha\cap|\tau|$, не имеющей самопересечений, и никакие $3$ точки среди $S$, $T$ и вершин ломаной $\alpha\cap|\tau|$ не лежат на одной прямой. Поэтому отрезок $ST$ пересекает $|\tau|\cap\left<\Delta\right>$ в нечетном количестве точек.
Аналогично прямая $\alpha\cap\left<\Delta\right>$ пересекает $|\tau|\cap\left<\Delta\right>$ в четном количестве точек.
Кроме того, никакие $4$ из вершин треугольника $\Delta$ и тетраэдра $\tau$ не лежат в одной плоскости.
Из всего этого получаем, что $\tau\cap\left<\Delta\right>:=\{\Gamma\cap\left<\Delta\right>\ :\ \Gamma\in\tau\}$ является \emph{$1$-циклом}, т.~е. набором отрезков в плоскости $\left<\Delta\right>$ таким, что каждая точка из $\left<\Delta\right>$ является концом четного (возможно, нулевого) количества отрезков набора.
В силу предположений об $\alpha$ мы можем выбрать отрезки 1-цикла так, чтобы прямая $\alpha\cap\left<\Delta\right>$ отделяла $O$ от $S,T$, и от всех вершин отрезков 1-цикла, отличных от $O$ (ср.~рис.~\ref{f:plane} слева).
Поэтому по аналогу леммы~\ref{0-eveng} для треугольника и $1$-цикла, некоторые два отрезка из $\partial\Delta$ и из $\tau\cap\left<\Delta\right>$ пересекаются, не имея общих вершин.
Не более одного из этих отрезков содержит точку $O$.
Поэтому найденный отрезок $1$-цикла $\tau\cap\left<\Delta\right>$ является пересечением с плоскостью $\left<\Delta\right>$ треугольника из $\tau$, пересекающего сторону треугольника $\Delta$, не имея с ней общих вершин.
\end{proof}

\subsection{Нереализуемость произведений в четырехмерном пространстве}\label{s:mnre4}

\emph{Зейфертовой цепью} (или кограницей) замкнутой ломаной $a$ в пространстве называется конечный набор $S$ треугольников (не вырожденных в отрезок или точку) в пространстве такой, что

$\bullet$ каждое ребро ломаной $a$ является стороной ровно одного треугольника из $S$;

$\bullet$ каждый отрезок, не являющийся ребром ломаной $a$, является стороной четного (возможно, нулевого) количества треугольников из $S$.

Две непересекающиеся замкнутые ломаные $a$ и $a'$ в пространстве \emph{зацеплены по модулю~$2$}, если для любых зейфертовых цепей $S$ для $a$ и $S'$ для $a'$
таких, что контур любого треугольника из $S$ не пересекается с контуром любого треугольника из $S'$, количество зацепленных по модулю~$2$ пар $(\Delta,\Delta')$ треугольников $\Delta$ из $S$ и $\Delta'$ из $S'$ нечетно.
Эквивалентность другим определениям зацепленности по модулю~$2$ (в частности, определению, приведенному перед теоремой~\ref{rlt-sachs}) доказана в~\cite[\S1.8]{Sk23'}, \cite[Утверждение 4.8.3]{Sk}.

\begin{proof}[Доказательство теоремы о квадрате \ref{0-ne4pr}]
Можно считать, что никакие $5$ из данных $25$ точек $A_{jp}$, $j,p\in[5]$, не лежат в одной трехмерной гиперплоскости (аналогично замечанию~\ref{r:stronger}.b).
Существует трехмерная гиперплоскость $\alpha$, по одну сторону от которой лежит ровно одна точка $A_{jp}$ из данных, а остальные 24 точки --- по другую (рис.~\ref{f:4spacepr}).
Точку пересечения гиперплоскости $\alpha$

$\bullet$ с отрезком $A_{jp}A_{kp}$ покрасим в синий для каждого $k\in[5]-\{j\}$;

$\bullet$ с отрезком  $A_{jp}A_{jq}$ покрасим в красный для каждого $q\in[5]-\{p\}$.

Пересечение произведения $jk\times pq$ с $\alpha$ назовем \emph{дугой}.
Тогда дуги имеют разноцветные концы.

\begin{figure}[h]\centering
 \begin{tikzpicture}[line cap=round,line join=round,>=triangle 45,x=0.32cm,y=0.32cm]
\clip(-4.93,-7.54) rectangle (15.43,4.42);
\fill[fill=black,fill opacity=0.1] (-3,-5) -- (0,1) -- (13,1) -- (10,-5) -- cycle;
\draw (-1,-4)-- (7,-2);
\draw (3.03,-1.9)-- (4.53,-2.43);
\draw (9,-4)-- (5.03,-2.61);
\draw (7,-2)-- (5.54,-1.51);
\draw (4,3)-- (1,0);
\draw (4,3)-- (-1,-4);
\draw (4,3)-- (3.03,-1.9);
\draw (4,3)-- (7,-2);
\draw (4,3)-- (9,-4);
\draw (4,3)-- (11,0);
\draw (2.39,-5)-- (2.03,-6.74);
\draw [dash pattern=on 2pt off 2pt] (3.03,-1.9)-- (2.39,-5);
\draw [dash pattern=on 2pt off 2pt] (7,-2)-- (8.65,-5);
\draw (8.65,-5)-- (9.47,-6.48);
\draw [dash pattern=on 2pt off 2pt] (1,0)-- (-2,-3);
\draw (-2,-3)-- (-3,-4);
\draw (-1.71,-5)-- (-2.4,-5.96);
\draw (9.71,-5)-- (10.61,-6.25);
\draw (12.24,-0.53)-- (13.38,-1.02);
\draw [dash pattern=on 2pt off 2pt] (-1.71,-5)-- (-1,-4);
\draw [dash pattern=on 2pt off 2pt] (9,-4)-- (9.71,-5);
\draw [dash pattern=on 2pt off 2pt] (11,0)-- (12.24,-0.53);
\draw (-3,-5)-- (0,1);
\draw (13,1)-- (10,-5);
\draw (10,-5)-- (-3,-5);
\draw (0,1)-- (1.33,1);
\draw (3,1)-- (3.34,1);
\draw (4,1)-- (5,1);
\draw (6,1)-- (8,1);
\draw (9.32,1)-- (13,1);
\draw (3.03,-1.9)-- (6.27,-1.13);
\draw (7.98,-0.72)-- (11,0);
\draw (1,0)-- (1.5,-0.17);
\draw (1.89,-0.3)-- (2.85,-0.62);
\draw (3.57,-0.86)-- (5.02,-1.34);
\draw (1,0)-- (-0.45,-0.72);
\draw (-0.45,-0.72)-- (-1,-4);
\draw (11,0)-- (10.01,-3.36);
\draw (10.01,-3.36)-- (9,-4);
\draw (4,3)-- (-0.45,-0.72);
\draw (4,3)-- (10.01,-3.36);
\draw [dash pattern=on 2pt off 2pt] (-0.45,-0.72)-- (-1.16,-1.32);
\draw [dash pattern=on 2pt off 2pt] (10.01,-3.36)-- (10.54,-3.92);
\draw (-1.16,-1.32)-- (-2.3,-2.27);
\draw (10.54,-3.92)-- (11.51,-4.96);
\draw [color=black] (1,0) ++(-1.5pt,0 pt) -- ++(1.5pt,1.5pt)--++(1.5pt,-1.5pt)--++(-1.5pt,-1.5pt)--++(-1.5pt,1.5pt);
\draw [color=black] (-1,-4) ++(-1.5pt,0 pt) -- ++(1.5pt,1.5pt)--++(1.5pt,-1.5pt)--++(-1.5pt,-1.5pt)--++(-1.5pt,1.5pt);
\fill [color=black] (7,-2) circle (1.5pt);
\draw [color=black] (9,-4) ++(-1.5pt,0 pt) -- ++(1.5pt,1.5pt)--++(1.5pt,-1.5pt)--++(-1.5pt,-1.5pt)--++(-1.5pt,1.5pt);
\draw [color=black] (11,0) ++(-1.5pt,0 pt) -- ++(1.5pt,1.5pt)--++(1.5pt,-1.5pt)--++(-1.5pt,-1.5pt)--++(-1.5pt,1.5pt);
\fill [color=black] (3.03,-1.9) circle (1.5pt);
\fill [color=black] (4,3) circle (1.5pt);
\draw[color=black] (3.05,3.72) node {$A_{11}$};
\fill [color=black] (2.03,-6.74) circle (1.5pt);
\draw[color=black] (3.45,-6.26) node {$A_{14}$};
\fill [color=black] (9.47,-6.48) circle (1.5pt);
\draw[color=black] (8.45,-6.34) node {$A_{13}$};
\fill [color=black] (-3,-4) circle (1.5pt);
\draw[color=black] (-3.91,-3.16) node {$A_{21}$};
\fill [color=black] (-2.4,-5.96) circle (1.5pt);
\draw[color=black] (-0.93,-5.8) node {$A_{31}$};
\fill [color=black] (10.61,-6.25) circle (1.5pt);
\draw[color=black] (12.35,-6.36) node {$A_{41}$};
\fill [color=black] (13.38,-1.02) circle (1.5pt);
\draw[color=black] (13.21,-2.14) node {$A_{51}$};
\fill [color=black] (-0.45,-0.72) circle (1.5pt);
\fill [color=black] (10.01,-3.36) circle (1.5pt);
\fill [color=black] (-2.3,-2.27) circle (1.5pt);
\draw[color=black] (-2.93,-0.94) node {$A_{12}$};
\fill [color=black] (11.51,-4.96) circle (1.5pt);
\draw[color=black] (13.41,-4.8) node {$A_{15}$};
\end{tikzpicture}
\caption{К доказательству теоремы о квадрате~\ref{0-ne4pr}}
\label{f:4spacepr}
\end{figure}

Аналогично последнему абзацу доказательства утверждения \ref{p:k4k4}.a либо

(1) в данном $(5,5)$-произведении найдутся два пересекающиеся треугольника, не имеющие общих вершин, либо

(2) любые две дуги пересекаются только в их общей вершине (если она есть).

В случае (2) пересечение плоскости $\alpha$ с телом данного $(5,5)$-произведения является PL изображением без самопересечений графа $K_{4,4}$ в $\alpha$.
Воспользуемся следующим PL аналогом теоремы~\ref{rlt-sachs}: \emph{в любом PL изображении без самопересечений графа $K_{4,4}$ в $\R^3$ найдутся два зацепленных по модулю 2 цикла длины 4}
(см. доказательство в \cite{Sa81, Z13}).
Получаем две зацепленные по модулю 2 замкнутые ломаные в $\alpha$, каждая из которых составлена из четырех дуг.
Возьмем  $\{a,b,a',b'\}=[5]-\{j\}$ и $\{c,d,c',d'\}=[5]-\{p\}$ такие, что дуги

$\bullet$ первой ломаной содержатся в произведениях $ja\times pc$, $jb\times pc$, $ja\times pd$, $jb\times pd$,

$\bullet$ второй ломаной содержатся в произведениях $ja'\times pc'$, $jb'\times pc'$, $ja'\times pd'$, $jb'\times pd'$.

Значит, ломаные являются пересечениями гиперплоскости $\alpha$ с телами $(3,3)$-произведений $jab\times pcd$ и $ja'b'\times pc'd'$.
Тогда получаем нужное по нижеприведенной лемме~\ref{0-even4g}.
\end{proof}

\begin{lemma}[о понижении размерности]\label{0-even4g}
Даны два $(3,3)$-произведения в $\R^4$, имеющие единственную общую вершину $O$.
Никакие $5$ из их $17$ вершин не лежат в одной трехмерной гиперплоскости.
Трехмерная гиперплоскость $\alpha$ отделяет $O$ от оставшихся $16$ точек $(3,3)$-произведений.
Эта гиперплоскость $\alpha$ пересекает их тела по двум непересекающимся замкнутым ломаным, зацепленным по модулю $2$ в $\alpha$.
Тогда в этих $(3,3)$-произведенияx найдутся два пересекающиеся треугольника, не имеющие общих вершин.
\end{lemma}

\begin{proof}
Обозначим через $S$ (через $S'$) множество всех тех треугольников из первого (из второго) $(3,3)$-произведений, которые не содержат точки $O$.
Так как никакие $5$ из $17$ вершин $(3,3)$-произведений не лежат в одной трехмерной гиперплоскости, то никакие $4$ из их $16$ проекций на $\alpha$ с центром $O$ не лежат в одной плоскости.
Поэтому контуры треугольников $\alpha\cap O\Delta$ и $\alpha\cap O\Delta'$ не пересекаются для любых $\Delta\in S$ и $\Delta'\in S'$.
Обозначим через $\gamma$ и $\gamma'$ непересекающиеся замкнутые ломаные из условия.
Так как $\gamma$ и $\gamma'$ зацеплены по модулю 2, то количество зацепленных по модулю 2 пар $(\Delta,\Delta')$ таких треугольников нечетно.
По лемме \ref{relation}.a такие треугольники зацеплены по модулю 2 тогда и только тогда, когда поверхности тетраэдров $O\Delta$ и $O\Delta'$ пересекаются в четном количестве точек (включая $O$).
Для любой стороны $MN$ треугольника из $S$, не лежащей в $\gamma$, и стороны $M'N'$ треугольника из $S'$, не лежащей в $\gamma'$, пересечение $OMN\cap OM'N'$ входит в ровно два пересечения боковых поверхностей  тетраэдров $O\Delta$ и $O\Delta'$.
Поэтому количество пар пересекающихся треугольников одного из трех следующих типов нечетно:

$\bullet$ пар $(\Delta,\Delta')$ для $\Delta\in S$ и $\Delta'\in S'$;

$\bullet$ пар $(OMN,\Delta')$ для стороны $MN$ ломаной $\gamma$ и треугольника $\Delta'\in S'$;

$\bullet$ пар $(\Delta,OM'N')$ для стороны $M'N'$ ломаной $\gamma'$ и треугольника $\Delta\in S$.

Теперь лемма следует из того, что в указанных парах треугольники не имеют общих вершин, а также содержатся в треугольниках из данных $(3,3)$-произведений.
\end{proof}

\emph{Книги, обзоры и учебные статьи в этом списке помечены звездочками.}